\documentclass[12pt]{amsart}
\usepackage{graphicx} 
\usepackage{pgfplots}
\pgfplotsset{compat=1.17}

\usepackage{amsmath,amsthm,amssymb,epsfig}
\usepackage{lmodern}
\usepackage{float}
\usepackage{url}
\usepackage{tikz}
\usepackage{mathtools}
\usepackage{subcaption}
\usepackage[linesnumbered,lined,boxed,commentsnumbered]{algorithm2e}
\usepackage{MnSymbol}
\usepackage{comment}
\usepackage{tkz-graph}
\usetikzlibrary{fit}
\usetikzlibrary{shapes}
\usepackage{enumitem}
\usepackage[left=3.5cm,top=3.5cm,right=3.5cm,bottom=3.5cm]{geometry}
\usepackage{mathrsfs}

\newtheorem{theorem}{Theorem}
\newtheorem{corollary}[theorem]{Corollary}
\newtheorem{lemma}[theorem]{Lemma}

\theoremstyle{definition}

\newcommand{\CL}{\mathrm{CL}}

\title{Between burning and cooling: liminal burning on graphs}

\author[A.\ Bonato]{Anthony Bonato}
\author[T.G.\ Marbach]{Trent G.\ Marbach}
\author[J.\ Marcoux]{John Marcoux}
\author[T.\ Mishura]{Teddy Mishura}
\address[A1,A2,A3,A4]{Toronto Metropolitan University, Toronto, Canada}

\email[A1]{(A1) abonato@torontomu.ca}
\email[A2]{(A2) trent.marbach@torontomu.ca}
\email[A3]{(A3) tmishura@torontomu.ca}
\email[A4]{(A4) jmarcoux@torontomu.ca}


\begin{document}
\input{tikz}
\begin{abstract} 
Liminal burning generalizes both the burning and cooling processes in graphs. In $k$-liminal burning, a Saboteur reveals $k$-sets of vertices in each round, with the goal of extending the length of the game, and the Arsonist must choose sources only within these sets, with the goal of ending the game as soon as possible. The result is a two-player game with the corresponding optimization parameter called the $k$-liminal burning number. For $k = |V(G)|$, liminal burning is identical to burning, and for $k = 1$, liminal burning is identical to cooling. 

Using a variant of Sperner sets, $k$-liminal burning numbers of hypercubes are studied along with bounds and exact values for various values of $k$. In particular, we determine the exact cooling number of the $n$-dimensional hypercube to be $n.$ 
We analyze liminal burning for several graph families, such as Cartesian grids and products, paths, and graphs whose vertex sets can be decomposed into many components of small diameter. We consider the complexity of liminal burning and show that liminal burning a graph is PSPACE-complete for $k\geq 2,$ using a reduction from $3$-QBF. We also prove, through a reduction from burning, that even in some cases when liminal burning is likely not PSPACE-complete, it is co-NP-hard. We finish with several open problems.
\end{abstract}

\maketitle
\section{Introduction}

Graph burning is a simplified model for the spread of influence or contagion in a network. The burning number of a graph, introduced in \cite{BJR0,BJR}, quantifies the speed at which the influence spreads to every vertex. Given a graph $G$, the burning process on $G$ is a discrete-time process that proceeds in rounds. Throughout the process, vertices are \emph{burned} or \emph{unburned}. At the beginning of the first round, all vertices are unburned. In each round, all unburned vertices that have a burned neighbor become burned, and then one new unburned vertex (if available) is chosen to burn; the latter vertices are called \emph{sources}. The \emph{burning number} of $G,$ written $b(G),$ is the minimum number of rounds needed until all vertices are burned.  Burning has since been extensively studied \cite{Burning_Hard,bessy2,bonf,bgs,bk,bl,bill,rj,kam,ko,LL,liu1,MitschePralatRoshanbin,prod}. For further background on graph burning, see the survey \cite{survey} and the book \cite{bbook}.

A recent variant of burning considers the dual problem, called \emph{cooling}, which maximizes the number of rounds for the influence to spread; see \cite{cooling}. Rather than vertices being burned or unburned, we refer to them as \emph{cooled} or \emph{uncooled}. Cooling is defined analogously as burning, except that the \emph{cooling number} of a graph $G$, written $\mathrm{CL}(G),$ is the maximum number of rounds needed until all vertices are cooled. In \cite{cooling}, the cooling number was studied for graph families such as paths, cycles, and Cartesian grids, as well as in models of social networks.

Given that burning and cooling are both extremal conditions, one may instead consider a middle ground. This was first considered in the \emph{burning game} in \cite{bill}, where two competing players make alternating selections of sources in each round. One player wants to minimize the number of rounds (like burning), while the other attempts to maximize the number of rounds (like cooling). In \cite{nir}, an unrelated adversarial burning game was studied.

We introduce a unified framework for both burning and cooling, generalizing both in a way distinct from the burning game. We define a two-player game called \emph{liminal burning} played on a graph $G$ with a positive integer parameter $k$. Play occurs over discrete rounds with opposing players, the \emph{Saboteur} and the \emph{Arsonist}. Vertices in $G$ occupy two binary states: \emph{unburned} or \emph{burned}, and \emph{revealed} or \emph{unrevealed}, where a vertex can be exactly one of unburned or burned, and exactly one of revealed or unrevealed. The game begins with all vertices of $G$ unrevealed and unburned. In each round, the Saboteur selects a set of $k$ unrevealed vertices and reveals them.  The Arsonist selects one of these $k$ revealed, unburned vertices to burn (if it exists). These are called \emph{moves} of the player. Burning spreads from any burned vertex to all of its neighbors, regardless of whether they are revealed or not; we refer to this spread as \emph{propagation}. If there are fewer than $k$ unburned and unrevealed vertices to choose from in a given round, the Saboteur must choose whichever vertices remain in those states. The play continues until all vertices are burned. 

\begin{figure}[htpb!]
    \centering
    \includegraphics[scale=0.4]{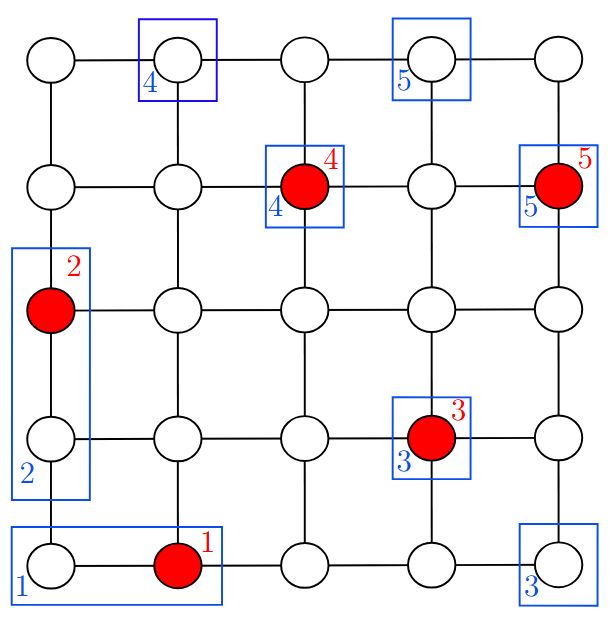}
    \caption{An example of a $2$-liminal burning game played on the $5 \times 5$ Cartesian grid. Sets are revealed by the Saboteur and marked in blue outlines, while vertices burned by the Arsonist are marked in red. The blue numbers refer to the order in which the sets are revealed, while the red numbers refer to the order in which the vertices are burned.}
    \label{fig1}
\end{figure}

For optimal play, the Saboteur extends the game as long as possible, while the Arsonist ends the game as soon as possible. We denote the number of rounds the process takes under optimal play by both players as $b_k(G)$, which we call the $k$-\emph{liminal burning number} of $G$.  See Figure~\ref{fig1} for an example of $k$-liminal burning with $k=2$.

Note that we may consider $k=f(n)$ to be constant or a function of the number of vertices $n$ of $G.$ A non-monotonic variant of liminal burning was first introduced in \cite{holden}, where revealed vertices become unrevealed in the next round. Unlike other pursuit-evasion games such as Cops and Robbers, the Arsonist always wins in the sense that all vertices eventually burn; however, the Saboteur heavily influences the value of $b_k$ through their choice of revealed vertices. 

As the name suggests, liminal burning occupies the space between burning and cooling. Observe that $b_1(G) = \CL(G)$, as the Saboteur will always choose vertices in a cooling sequence to reveal in each round. Also, $b_{|V(G)|}(G) = b(G),$ as all vertices are revealed in the first round. The following inequalities are immediate:
$$b(G) = b_{|V(G)|}(G) \leq b_{|V(G)|-1}(G) \leq \cdots \leq b_2(G) \leq b_1(G) = \CL(G).$$

The paper's main results are as follows. In Section~2, we study liminal burning in hypercubes. Using variations of Sperner sets in finite set theory, we prove that the cooling number of the hypercube $Q_n$ is $n$ in Theorem~\ref{thm:CLQN}, answering an open problem in \cite{cooling}. We extend these results for values $k\ge 2$, and establish lower and upper bounds for $b_k(Q_n)$ in Corollary~\ref{thm:bkqn} and Theorem~\ref{thm:up1}. In Section~3, we consider bounds on the $k$-liminal burning number for general graphs. A lower bound for $b_k$ in terms of the cooling number is given in Theorem~\ref{LB}. We introduce $(k,d)$-special graphs and use these to give upper bounds for $b_k$ on caterpillar graphs. We also give an upper bound in terms of the edge cover of an associated hypergraph of the graph in Theorem~\ref{thm:edge_cover}. In Section~4, we establish bounds for $b_k$ on Cartesian products of graphs. We determine $b_k$ for Cartesian grids within a constant ratio for all values of $k$ (see Theorems~\ref{gridLB} and \ref{grub}) and summarize results on $b_k$ of Cartesian grids in Table~1. In Section~5, we consider the complexity of computing $b_k$. We show that for all $k\ge 2,$ computing $b_k$ is PSPACE-complete using a reduction from $3$-QBF in Theorem~\ref{thm:pspace}. In some cases where liminal burning is likely not PSPACE-complete, it is shown to be co-NP-hard through a reduction from burning. The final section considers open problems.

We consider only simple, finite, and undirected graphs. In a graph $G$, we denote the distance between vertices $u$ and $v$ by $d_G(u,v)$, writing $d(u,v)$ if $G$ is clear from context. The diameter of $G$ is $\mathrm{diam}(G)$. The complete graph of order $n$ is $K_n$. For a positive integer $n$, we let $[n] = \{ 1, 2, \ldots ,n\}$. The power set of a set $X$ is $\mathcal{P}(X)$. The \emph{hypercube} with dimension $n,$ written $Q_n,$ has vertices $\mathcal{P}([n])$, with vertices $u$ and $v$ adjacent if and only if $|u \Delta v| =1$, where $\Delta$ is the symmetric difference. Note that in hypercubes, $d(u,v) = |u \Delta v|$. 
For background on graph theory, see \cite{west}. 

\section{Hypercubes}

We consider liminal burning in hypercubes. In Theorem~\ref{thm:CLQN}, we determine the precise value of $\mathrm{CL}(Q_n)$, which was an open problem in \cite{cooling}. To do so, we use certain extremal families of sets that will be useful for later results on $b_k(Q_n)$ in this section.

A \emph{Sperner family} is a collection of subsets $\mathcal{F}$ of a finite set $\Omega$ such that none of the subsets contains another. We consider a variant of this where each set has a distinct cardinality, not including the sets of cardinality $0$ and $|\Omega|$, as these sets preclude the inclusion of any other.
We also do not include any sets of cardinality $1$, as if $\{d\} \in \mathcal{F}$, then no other subset in the family may contain $d$, and so this is a trivial case. 

A family of sets $\mathcal{F}$ that are subsets of $[n]$ is a $k$-\textit{rainbow Sperner family on $[n]$ of cardinality} $j$ if the following conditions are satisfied:
\begin{enumerate}
\item $|\{S \in \mathcal{F}:|S| = i\}| = k$ for all $2 \leq i \leq j + 1$.
\item For all $S,T \in \mathcal{F}$ with $S \neq T$, $S \not\subseteq T$. 
\end{enumerate}
If $k=1$, then we say that $\mathcal{F}$ is a \textit{rainbow Sperner family on $[n]$ of cardinality} $j$. Rainbow Sperner families always exist, as the next result shows.

\begin{lemma}\label{lemmr}
Given $n\geq 3$, there exists a rainbow Sperner family on $[n]$ of cardinality $n-2$.     
\end{lemma}
\begin{proof}
We will prove this result using induction on $n$. The two base cases are for $n=3$ with $\{1,2\}$, and $n=4$ with $\{1,2\},\{1,3,4\}$. For the sake of induction, suppose that for $[n-2]$, there exist $A_2, A_3, \ldots, A_{n-3}$ that satisfy the conditions of the lemma. Define $B_2 = \{n-1,n\}$, $B_{i+1}=A_i \cup \{n\}$ for $2 \leq i \leq n-4$, and $B_{n-1}=[n] \setminus \{n\}$. Each $B_i$ is contained in $[n]$, and each has cardinality $|B_i|=i$. 

We now show that $B_i \not\subseteq B_j$ for each $i \neq j$. 
If $B_i \subseteq B_j$ for $3 \leq i,j \leq n-2$ and $i \neq j$, then since $n \in B_i \cap B_j$ and $n-1$ is in neither, it must be that $A_{i-1} \subseteq A_{j-1},$ contradicting the definition of the $\{A_i\}$. 
The subset $B_{2}$ is the only subset that contains both $n-1$ and $n$, so it cannot be contained in any $B_i$, and trivially cannot contain any $B_i$. 
The subset $B_{n-1}$ likewise cannot be contained in $B_i$. 
If some $B_i\subseteq B_{n-1}$, then $B_i$ could not contain $n$, since $n \notin B_{n-1}$.
However, all other subsets contain $n$, so no such $B_i$ exists, so none of the $B_i$ contains another. 
\end{proof}

The rainbow Sperner family in Lemma~\ref{lemmr} is of the maximum cardinality possible for all $n\geq 3$. For $n=2$, we can include one of the two trivial subsets, say $\{1,2\}$, to find a family of cardinality $n-1=1$. Additionally, we present the following lemma about the relation in symmetric difference for rainbow Sperner families. 

\begin{lemma}\label{lem:rainbow_sperner_symdif}
    Let $\mathcal{F} = \{S_2,S_3,\ldots,S_{n-1}\}$ be a rainbow Sperner family on $[n]$ of cardinality $n-2$. For all $2 \leq i\neq j \leq n-1,$ we have that $|S_i \Delta S_j| \geq |i-j|+2$. 
\end{lemma}

\begin{proof}
    Consider two sets $S_i,S_j \in \mathcal{F}$. As neither set is a proper subset of the other, they share at most $\min(i,j)-1$ elements in common. Thus, we have that 
    \begin{eqnarray*}
|S_i \Delta S_j| &\geq& (\max(i,j) -\min(i,j) +1) + (\min(i,j)-\min(i,j)+1) \\
&=& |i-j|+2.
\end{eqnarray*}
The proof follows.
\end{proof}

Observe that $\CL(Q_1) = 2$. The remaining cases of $n>1$ are covered in the following theorem, which, for the first time, gives the cooling number of the hypercube. Our proof heavily relies on the existence of rainbow Sperner families. 

\begin{theorem}\label{thm:CLQN}
If $n\geq 2$, then $b_1(Q_n)=\CL(Q_n) = n.$
\end{theorem}

\begin{proof}
For $n=2$, a cooling sequence of $(0,0),(1,1)$ yields that $\CL(Q_2)=2$, and so we may assume that $n \geq 3$.

We first claim that $\CL(Q_n)\leq n$. 
Let $u$ be the first cooled vertex in the cooling process, and let $v$ be the unique vertex of distance $n$ from $u$.
All vertices with symmetric difference at most $i$ from $u$ will be cooled after propagation in round $i+1$ due to cooling $u$ in the first round.
Thus, $v$ is the only vertex that may not be cooled after propagation in round $n$, so $v$ must be selected and cooled at the end of round $n$. 
Therefore, all vertices are cooled by the end of round $n$, from which it follows that $\CL(Q_n)\leq n$. 

We show that $\CL(Q_n) \geq n$, which completes the proof. By Lemma~\ref{lemmr}, there exists a rainbow Sperner family $S_2,\ldots,S_{n-1}$ of cardinality $n-2$ on $[n]$. We claim that $(\emptyset,S_2,\ldots,S_{n-1})$ is a cooling sequence of $Q_n$. By Lemma~\ref{lem:rainbow_sperner_symdif}, we know that $|S_i \Delta S_j| \geq |i-j|+2$ for $2 \leq i\neq j \leq n-1$, and as $|\emptyset \Delta S_i| = i = |i-1|+1,$ we find that this is indeed a cooling sequence. Furthermore, we note that $[n]$ must be unburned at the end of round $n-1$, as $|[n] \Delta S_i| = n-i$ for $2 \leq i \leq n-1$, and so this sequence results in $Q_n$ being burned in exactly $n$ rounds, as claimed.
\end{proof}
For $b_k(Q_n)$ with general values of $k,$ we derive the following lower bound that is an application of Theorem~\ref{thm:CLQN}.
\begin{corollary}\label{thm:bkqn}
If $k$ is a positive integer and $Q_n$ is the graph of the $n$-dimensional hypercube, then $b_k(Q_n) \geq n - \left\lceil \log_2 k \right\rceil.$ 
\end{corollary}
\begin{proof}
Consider the graph $Q_{n-\left\lceil \log_2 k\right\rceil}$, and let $K = [n]\setminus\left[n-\left\lceil \log_2 k \right\rceil\right]$. By Theorem~\ref{thm:CLQN}, we know that there exists a cooling sequence $c_1,c_2,\ldots,c_{n-\left\lceil \log_2 k\right\rceil-1},$ where each $c_i \subseteq \left[n-\left\lceil \log_2 k \right\rceil\right]$. We propose the following strategy for the Saboteur: in round $i$, reveal $k$ arbitrarily chosen vertices $V_i \subseteq c_i \cup \mathcal{P}(K)$. We note that $d(V_i,V_j) \geq |c_i \Delta c_j| \geq |i-j|+1$, so these vertices can always be revealed no matter the choice of the Arsonist. Additionally, in round $n-\left\lceil \log_2 k\right\rceil-1$, the vertex $[n]$ is unburned. Thus, the game lasts at least one more round, which gives the desired lower bound. 
\end{proof}
The bound in Corollary~\ref{thm:bkqn} is likely far from tight for general $k$. We next determine values of $b_k(Q_n)$ precisely for $k=2,3,4$.
\begin{theorem}\label{thm:b2qn}
    If $n \geq 7$, then $b_2(Q_n)=n-1$.
\end{theorem}
\begin{proof}
Note that by Corollary~\ref{thm:bkqn}, $b_2(Q_n) \geq n-1$ for $n \geq 7$, so we need only show that $b_2(Q_n) < n$ for $n \geq 7$. For a contradiction, suppose that $b_2(Q_n) = n$ and consider an optimal $2$-liminal burning game played on $Q_n$. There must exist sets $S_1 = \{v_1,v_{1}'\},\ldots,S_{n-1}=\{v_{n-1},v_{n-1}'\}$ revealed by the Saboteur and vertices $a_1,a_2,\ldots,a_{n-1}$ burned by the Arsonist, with $a_i \not \subseteq v_j, v_j'$ for $j > i$. We assume without loss of generality that $a_1 = \emptyset$, as the Saboteur can reveal $\emptyset$ and another vertex in the first round, observe the Arsonist's choice, then relabel the graph. 
    
Consider now $v_{n-1}$ and $v_{n-1}'$. We know that $|v_{n-1}|,|v_{n-1}'|\geq n-1$ and without loss of generality let $v_{n-1}$ be a set of $n-1$ elements such that $1 \notin v_{n-1}$. This implies that $1 \in a_i$ for $2 \leq i \leq n-2$, and, letting $a_2 = \{1,y\}$, we have that $y \notin v_{n-2}, v_{n-2}'$. Thus, $v_{n-2} = [n]\setminus \{y,z_1\}$ and $v_{n-2}'=[n]\setminus\{y,z_2\}$, where $z_1 \neq z_2 \neq 1$. However, for both $v_{n-2}$ and $v_{n-2}'$ to have been revealed by the Saboteur, we know that $a_3,a_4 \not \subseteq v_{n-2},v_{n-2}'$. Noting that $y \notin a_3,a_4$, we derive that $z_1,z_2 \in a_3,a_4$, which forces $a_3 = \{1,z_1,z_2\}$ because $1 \in a_3,a_4$. Finally, we find that $a_3 \subseteq a_4$, which is a contradiction, completing the proof.
\end{proof}

\begin{theorem}\label{thm:b4qn}
If $n\ge 11,$ then $b_4(Q_n) = n-1$. 
\end{theorem}

\begin{proof}
Let $X_2,X_3,\ldots,X_{n-8}$ be a rainbow Sperner family on $[n-4]\setminus\{1,2,3\}$ where each $X_i$ contains 4 except for $X_{n-8}$, let $Y_2,\ldots,Y_{n-8}$ be a rainbow Sperner family on $[n-4]\setminus\{1,2,4\}$ where each $Y_i$ contains 3 except for $Y_{n-8}$, let $Z_2,\ldots,Z_{n-8}$ be a rainbow Sperner family on $[n-4]\setminus\{1,3,4\}$ where each $Z_i$ contains 2 except for $Z_{n-8}$, and finally let $W_2,W_3,\ldots,W_{n-8}$ be a rainbow Sperner family on $[n-4]\setminus\{2,3,4\}$ where each $W_i$ contains 1 except for $W_{n-8}$. The strategy for the Saboteur is as follows: In round 1, reveal the vertices $\emptyset, \{1\}, \{2\},$ and $\{3\}$. By relabeling if necessary, we can assume without loss of generality that the Arsonist burns vertex $\emptyset$ and that at the beginning of round 2, the vertices $\{1,2\}$ and $\{1,3\}$ are revealed. For rounds $2 \leq i \leq 4$, the Saboteur reveals the set of vertices $S_i$, where
\begin{align*}
S_2 &=\left\{\{1,4\}, \{1,n-3\}, \{2,3\}, \{2,4\}\right\}, \\
S_3 &=\left\{\{2,n-3,n-2\}, \{2,n-2,n-1\}, \{2,n-2,n\}, \{3,n-3,n-1\}\right\},\text{ and}\\
S_4 &=\left\{\{3,n-2,n-1,n\}, \{4,n-3,n-2,n\}, \{4,n-3,n-1, n\}, \{4,n-2,n-1,n\}\right\}
\end{align*}
For rounds $5 \leq i \leq n-5$ and round $n-4$, the Saboteur reveals the set of vertices $S_i$, where
\begin{equation*}
S_i = \begin{cases}
X_{i-3}\cup\{n-3,n-2,n-1\} \\ 
Y_{i-3}\cup\{n-3,n-2,n\} \\
Z_{i-3}\cup\{n-3,n-1,n\} \\ 
W_{i-3}\cup\{n-2,n-1,n\}         
\end{cases}
\text{and}\quad
S_{n-4} = \begin{cases}
[n]\setminus\{1,2,n-3,n-2\} \\ 
[n]\setminus\{1,2,n-3,n-1\} \\ 
[n]\setminus\{1,2,n-3,n\} \\ 
[n]\setminus\{1,2,n-2,n-1\}.
\end{cases}
\end{equation*}
We note that the $S_i$ form a $4$-rainbow Sperner family on $[n]$ of cardinality $n-4$, and note as well that $\{1,2\}$ and $\{1,3\}$ are not proper subsets of any element of $S_i$, so the strategy is possible. Furthermore, as only vertices of cardinality at most $n-4$ have been burned by round $n-4$, no vertices of cardinality $n-2$ are burned at the start of round $n-3$. Thus, the Saboteur can reveal any four vertices $w_1,\ldots,w_4$ of cardinality $n-2$ in round $n-3$. Additionally, as $n \geq 11$, there exists $v \in V(Q_n)$ of cardinality $n-1$ that does not contain $w_1,w_2,\ldots,w_4$. Thus, $v$ and $[n]$ are unburned at the start of round $n-2$, so the Saboteur may reveal those in round $n-2$. At least one more round is needed to burn the other revealed vertices, so $b_4(Q_n) = n-1$. 
\end{proof}

The following corollary is immediate.

\begin{corollary}
If $n\ge 11,$ then $b_3(Q_n)=n-1$. 
\end{corollary}

Additionally, we suspect that if $n$ is large enough compared to $k$, this trend continues. It seems as though once $n$ is large enough, there will be enough ``space'' to construct $k$-rainbow Sperner families of necessary cardinality. We conjecture that if $k$ is a positive integer, then there exists an integer $n_k$ such that $b_k(Q_n) = n-1$ for all $ n \geq n_k$. 

By Theorem~\ref{thm:CLQN}, for all $k\ge 1,$ we have that $b_k(Q_n)\leq n$. We can improve this bound using results of Alon~\cite{alon} and Kleitman~\cite{Kleitman}. 
\begin{theorem}[\cite{alon}]
\label{thm:alon_burning_result}
The burning number of the hypercube $Q_n$ equals $\left\lceil\frac{n}{2}\right\rceil +1$, which may be achieved by choosing any vertex $u$ as a source in the first round and the unique vertex of maximum distance from $u$ as a source in the second round.  
\end{theorem}
\begin{theorem}[\cite{Kleitman}] \label{thm:Kleitman}

Let $d,n$ be natural numbers such that $n \geq 2d+1$, and suppose we have a $F\subseteq V(Q_n)$ that satisfies $\max_{u,v \in F}d(u,v) \leq 2d$. We then have that
\[ |F| \leq \sum _{i=0}^{d} \binom{n}{i}.\]
\end{theorem}
We have the following upper bound.
\begin{theorem}\label{thm:up1}
Suppose $n+2 \leq k < 2^{n-1}$ and suppose $d$ is the largest integer such that $k > \sum_{i=0}^{d}\binom{n}{i}$. We then have that
\[b_k(Q_n) \leq n-d+1.\]
\end{theorem}
\begin{proof}
We force the following two conditions onto the Arsonist, noting that both only benefit the Saboteur. First, the Arsonist must choose the first two vertices they will burn during the first round, and they must be of maximum possible distance within the set of $k$ vertices revealed by the Saboteur. Second, the Arsonist cannot burn other vertices for the rest of the game.

Without loss of generality, say the first vertex is $\emptyset$ and the second is $u$. 
By Theorem~\ref{thm:Kleitman}, as $u$ and $\emptyset$ are maximally distant, $|u\Delta\emptyset| > 2d$.
Since $k\geq n+2$, there must be a valid choice for $u$ for the Arsonist's second move.

The vertices on a shortest path between $\emptyset$ and $u$ form a subhypercube of dimension at least $2d+1$. 
The vertices of this subhypercube will be burned by the end of round $\left\lceil \frac{2d+1}{2}\right\rceil+1 = d+2$ by Theorem~\ref{thm:alon_burning_result}. 
Each remaining vertex in the graph has distance at most $n-(2d+1)$ from some vertex on the subhypercube, and so each of these will be burned by the end of round $(d+2) + n-(2d+1) = n-d+1$. 
\end{proof}

\section{Bounds}

We next examine some general bounds on the $k$-liminal burning number. Additionally, we investigate connections between the burning and cooling numbers and the $k$-liminal burning number. Observe that after $\left\lceil \frac{n}{k} \right\rceil $ rounds in the game, all vertices are available to the Arsonist. This gives an immediate bound.

\begin{theorem} \label{lem:nk+b}
For a graph $G$ on $n$ vertices, $b_k(G) \leq b(G) + \left \lceil \frac{n}{k} \right\rceil$.
\end{theorem}
For example, Theorem~\ref{lem:nk+b} gives that if $k = \omega\left(\frac{n}{b(G)}\right)$, then $b_k(G) = (1 + o(1))b(G)$. 
A less obvious bound connects the cooling number and the $k$-liminal burning number.

\begin{theorem}\label{LB}
For a graph $G$ and all positive integers $k$, we have that
\[\left\lceil\frac{\CL(G)}{k}\right\rceil \leq b_k(G).\]
\end{theorem}

If $k$ is a constant, then Theorem~\ref{LB} gives that $\left\lceil\frac{\CL(G)}{k}\right\rceil \le b_k(G) \le \CL(G)$, which shows that $b_k(G)$ is within a constant multiple of $\CL(G)$.  

\begin{proof}[Proof of Theorem~\ref{LB}]
We show that the Saboteur has a strategy that lasts $\left\lceil\frac{\CL(G)}{k}\right\rceil$ rounds independent of the Arsonist. Given an optimal cooling sequence $v_1,v_2,\dots$ $v_{t}$, the Saboteur reveals $v_1,v_2,\dots,v_k$, then $v_{k+1},v_{k+2},\dots,v_{2k}$, and so on until the entire cooling sequence has been revealed. If $\CL(G)$ is not a multiple of $k$, then the Saboteur can reveal an arbitrary set of extra vertices to fill out their selection. We claim that this strategy of the Saboteur is always possible. 

In the first case, suppose there is a cooling sequence of length $\CL(G)$. Therefore, for all $1 \leq i < j \leq t$, we have that $d(u_i,u_j) \geq j - i + 1$; otherwise, in round $j$, vertex $u_j$ has already been cooled by the spread of $u_i$. Suppose the Saboteur has so far revealed $v_1,v_2,\dots,v_{rk}$ for some $1 \leq r \leq \left\lceil \frac{\CL(G)}{k}\right\rceil - 1$, and it is the Saboteur's move. If $v_{rk+1},\dots,v_{(r+1)k}$ is not a valid choice, then there is some $1 \leq t \leq r - 1$ with  $tk + 1 \leq i \leq (t+1)k \leq rk$ and $rk + 1 \leq j \leq (r+1)k$ with $d(v_i,v_j) < r - t$. Hence, $v_i$ was revealed on or following round $t$ and has spread to cover $v_j$ by the beginning of round $r$. Since the $v_i$'s form a cooling sequence and $k \geq 1$, we have that 
\begin{equation*}
    d(v_i,v_j) \geq j - i + 1 \geq rk + 1 - tk - 1 + 1 = rk - tk + 1 > rk-tk \geq r-t.
\end{equation*}
As this contradicts our assumption, the result follows in this case.

In the second case, suppose there are no cooling sequences of length $\CL(G)$, so the longest cooling sequence is of length $\CL(G) - 1$. In this case, suppose we give the Arsonist the advantage that they can burn $k$ vertices every round and keep the strategy of the Saboteur the same. We need only consider the case where $\CL(G) - 1$ is a multiple of $k$. After $\frac{\CL(G) - 1}{k}$ rounds, we have that every vertex in the cooling sequence is burned, and they burned for at most as long as if we were cooling them in sequence. In that case, the graph is not fully cooled until the next round, so the graph cannot be fully burned here. Thus, the game lasts at least $\frac{\CL(G)-1}{k} + 1 =\left\lceil \frac{\CL(G)}{k}\right\rceil $ rounds.
\end{proof}

Note that Theorem~\ref{LB} is unlikely to be tight for larger $k$ since any pair of vertices cooled $t$ rounds apart must be of distance at least $t+1$ from each other, so these $t$ intermediate vertices can be revealed by the Saboteur to slow the process further. We next give bounds on $b_k$ for paths $P_n$ with $n$ vertices and show that they are asymptotically tight for certain instances of $k$ in Corollary~\ref{pathc}.

\begin{theorem} \label{thm:paths}
For positive integers $k$ and $n$, we have that \[\left\lfloor \frac{n}{k+1} \right\rfloor + \left\lfloor\frac{-1 + \sqrt{5+4k}}{2} \right\rfloor \leq b_k(P_n) \leq \left\lceil \frac{n}{k} \right\rceil + k - 1. \] 
\end{theorem}

\begin{proof}
Let $V(P_n)$ be $\{v_1,v_2,\dots,v_n\}$ with $v_i$ adjacent to $v_{i+1}$ for $1 \leq i \leq n-1$.

For the lower bound, consider the strategy where the Saboteur always reveals the $k$ unrevealed and unburned vertices which are closest to $v_1$. To ease our analysis, we will only require the Arsonist to burn the last $k$ vertices of the path. As a consequence of this advantage, the Arsonist always chooses the largest available index, since if they choose anything else, it will take longer for there to be any vertices among the last $k$ available to be burned, which only benefits the Saboteur. Thus the Saboteur will reveal $v_{(k+1)(t-1) + 1},v_{(k+1)(t-1) + 2},\dots,v_{(k+1)(t-1) + k}$ in round $t$ and the Arsonist will burn $v_{(k+1)(t-1) + k}$. The Saboteur continues like this until there are fewer than $k + 1$ vertices to reveal, and this leaves at least $k$ vertices at the end of the path that need to be burned. This strategy takes $\left\lfloor \frac{n}{k+1}\right\rfloor$ rounds, giving the first part of the bound. 
    
For simplicity, we will assume these remaining vertices are the last $k$ vertices on the path. To burn the last $k$ vertices, the Arsonist takes advantage of the fact that the vertex $v_{n-k}$ is already burned, so we can allow burning to spread and find the optimal number of rounds to burn these vertices given this help. This is equivalent to finding the smallest positive $t$ satisfying 
\[
t - 1 + \sum_{i=0}^{t-1}(2i+1) \geq k,
\]
which is $t = \frac{-1 + \sqrt{5 + 4k}}{2}$. Combining these two parts of the strategy gives the desired lower bound.
    
To prove the upper bound, we partition the vertex set of $P_n$ into the $k$-sets $S_i = \{v_{(i-1)k+1},\dots,v_{ik}\}$ for $1 \leq i \leq \left\lceil \frac{n}{k} \right\rceil$, letting $S_{\left\lceil \frac{n}{k}\right\rceil}$ be smaller when necessary. We claim that in every round, the Arsonist can place a source in some $S_i$ that does not yet contain a source. First, this is certainly possible in the first round, so we proceed by induction. In the $(t+1)^{th}$ round, the Saboteur has revealed a total of $k(t+1)$ vertices, and the Arsonist has chosen sources in $t$ unique sets. As long as the Saboteur is forced to reveal a vertex outside those sets, the Arsonist can continue. The $t$ sets contain a total of at most $kt$ vertices, and the Saboteur has so far revealed $kt + k$ vertices, so at least one of them must be outside those $t$ sets, and so the Arsonist has a valid choice. Hence, the Arsonist does this until all the sets contain a source, then waits $k-1$ rounds, and all vertices are burned.
\end{proof}

Given Theorem~\ref{thm:paths}, we have the following bounds on $b_k(P_n)$ for different values $k$ in terms of $n$.
\begin{corollary}\label{pathc}
Let $n$ be a positive integer and let $k$ be a positive integer-valued function of $n$.
\begin{enumerate}
    \item If $k = O(1)$, then $b_k(P_n) = \Theta(n).$
    \item If $k = \omega(1) \cap o(\sqrt{n})$, then $b_k(P_n) = (1 + o(1))\frac{n}{k}.$
    \item If $k = \Theta(\sqrt{n})$, then $b_k(P_n) = \Theta(\sqrt{n}).$
    \item If $k = \omega(\sqrt{n}) \cap o(n)$, then $b_k(P_n) = (1+o(1))\sqrt{n}.$
    \item If $k = \Theta(n)$, then $b_k(P_n) = \sqrt{n} + O(1).$
\end{enumerate}
\end{corollary}

\begin{proof}
    Item (1) follows from Theorem~\ref{LB}. Items (2) and (3) follow from Theorem~\ref{thm:paths}. Items (4) and (5) follow from Theorem~\ref{lem:nk+b}.
\end{proof}

Interestingly, we did not determine asymptotically tight values of $b_k(P_n)$ for all ranges of $k$. A step towards settling this would be to settle whether
$$b_k(P_n) \leq \left\lceil \frac{n}{k+1} \right\rceil + \left\lceil \frac{k-1}{2} \right\rceil.$$

We also note the following corollary that establishes that non-consecutive liminal burning parameters can be arbitrarily far apart.
\begin{corollary}
For all positive integers $k$, $t$ and $j$ with $j \geq 2$, there exists a graph $G$ such that $$|b_{k}(G) - b_{k + j}(G)| > t.$$
\end{corollary}

\subsection{$(k,d)$-special graphs}
We provide conditions that give similar upper bounds to the one in Theorem~\ref{thm:paths} for paths, but for more general classes of graphs. A graph $G$ is $(k,d)$-\emph{special} if its vertices can be partitioned into sets $S = \left\{S_1,S_2,\ldots,S_{\left\lceil \frac{n}{k}\right\rceil}\right\}$ of cardinality $k$ (with at most one exception) such that all sets have diameter at most $d$ with respect to distance in $G$. We call such a set $S$ a \emph{witness} of $G$ being $(k,d)$-special. If $k = d$, then we say $G$ is $k$-\emph{special}. The following theorem gives an upper bound on $b_k(G)$ if $G$ is $(k,d)$-special.

\begin{theorem} \label{thm:kd_special_upperbound}
If $G$ is $(k,d)$-special and order $n$, then $b_k(G) \leq \left\lceil \frac{n}{k} \right\rceil + d.$
\end{theorem}

\begin{proof}
Let $S$ be a witness of $G$ being $(k,d)$-special. We claim that in every round, the Arsonist can place a source in some set in $S$ that does not yet contain a source. This is possible in the first round, so we proceed by induction on the number of rounds. In the $(t+1)^{th}$ round, the Saboteur has revealed $k(t+1)$ vertices, and the Arsonist has chosen sources in $t$ unique sets from $S$. As long as the Saboteur must reveal a vertex outside those sets, the Arsonist can continue. The $t$ sets contain at most $kt$ vertices, and the Saboteur has so far revealed $kt + k$ vertices. Hence, at least one of them must be outside those $t$ sets, so the Arsonist does this until all the sets in $S$ contain a source, then waits $d$ rounds for all vertices to burn.
\end{proof}

Note that not all graphs are $k$-special for every value of $k$. A \emph{spider} is a tree with exactly one vertex of degree larger than two; such a vertex is called the \emph{head}. For $k \geq 3$, consider a spider $G$ with $k + 1$ paths connected to the head each having $k-1$ vertices. $G$ has $k^2$ vertices, so there must be $k$ sets in a witness of it being $k$-special. However, there are $k+1$ leaves all at a distance of $2k$ from each other, so they must all be in different sets, a contradiction. For $k=2$, we have the following.

\begin{theorem}\label{ttree}
If $G$ is a tree of order $n$, then it is $2$-special, and so satisfies $b_2(G) \le \left\lceil \frac{n}{2} \right\rceil + 2.$
\end{theorem}

\begin{proof}
We note that $K_1$ and $K_2$ are $2$-special and proceed by induction. 
Let $u$ be a vertex that is adjacent to at least one leaf and adjacent to exactly one non-leaf. 
Note that such a vertex always exists: Consider any vertex that would be a leaf after all leaves of $G$ are deleted. 
If $u$ is adjacent to an even number of leaves, then pair up the leaves, and if there is an odd number of leaves, then pair up all but one and pair the last one with $u$. We can ignore these vertices and proceed with the rest of the tree.
\end{proof}

Note that while Theorem~\ref{ttree} implies that the structure of $k$-special graphs is only interesting when $k\ge 3$, the implied upper bound is similar to one of the upper bounds given on the cooling number in~\cite{cooling} and so is unlikely to be tight.

We can apply Theorem~\ref{thm:kd_special_upperbound} to caterpillars. A \emph{caterpillar} is a tree that becomes a path when all leaves are removed; that path is called the \emph{spine} of the caterpillar. 
    
\begin{theorem}
For all positive integers $k$, caterpillars are $k$-special. Hence, a caterpillar $G$ of order $n$ satisfies $b_k(G) \le \left\lceil \frac{n}{k} \right\rceil + k.$
\end{theorem}

\begin{proof}
If $k=1$, then this is immediate. If $k=2$, then it is true because all caterpillars are trees, so suppose $k \geq 3$.
        
Let $v_1,v_2,\dots,v_r$ be the spine of the caterpillar, and let $u^i_j$ be the leaves adjacent to $v_i$ for $1 \leq j \leq \mathrm{deg}(v_i)$. Order the vertices as $$u^1_1,\dots,u^1_{\mathrm{deg}(v_1) - 1},v_1,u^2_1,\dots,u^2_{\mathrm{deg}(v_2) - 2},v_2,\dots,u^r_{\mathrm{deg}(v_r)-1},v_r.$$ We refer to these vertices as $w_1,w_2,\dots,w_n$. We claim that the sets $$S_i = \{w_{(i-1)k+1},\dots,w_{ik}\}$$ for $1 \leq i\leq \left\lceil \frac{n}{k} \right\rceil$ with $S_{\left\lceil \frac{n}{k} \right\rceil}$ being smaller if necessary, gives the desired decomposition. 

Let $i$ be a nonnegative integer, and consider $S_i$. First suppose $w_{(i-1)k+1}$ is on the spine. If $S_i$ induces a connected subgraph, then we are done. Hence, we can assume there is some nonnegative $t$ such that $w_{ik-t}$ and $w_{ik}$ are both leaves, are both in $S_i$, and share a neighbour which is not $w_{(i-1)k+1}$. All of these vertices have the same distance to the rest of $S_i$. We can assume, therefore, without loss of generality, that there is no positive $t$ such that $w_{ik-t}$ and $w_{ik}$ are both leaves and share a neighbor.

Now observe that the distance between $w_{(i-1)k+1}$ and $w_{ik}$ is at most $k$ since the distance from $w_{(i-1)k+1}$ to $w_{ik - 1}$ is at most $k-2$. This follows since $S_i \setminus \{w_{ik}\}$ induces a connected subgraph and the distance from $w_{ik - 1}$ to $w_{ik}$ is $2$. 
If $w_{(i-1)k+1}$ is a leaf, then either all vertices in $S_i$ are leaves adjacent to the same spine vertex, or there is some vertex in $S_i \setminus \{w_{(i-1)k + 1}\}$ which is a spine vertex. In the first case, all the vertices are a distance of $2$ from each other. In the second case, we can use the next spine vertex in the ordering and apply a similar argument to the case where $w_{(i-1)k+1}$ was a spine vertex. 
\end{proof}
Being $(k,d)$-special can be viewed in terms of hypergraph matchings. For a graph $G$, let $H_{k,d}(G)$ be the hypergraph such that $ V(H_{k,d}(G))=V(G)$, and edges consisting of all subsets of $V(G)$ of cardinality at most $k$ and diameter at most $d$ with respect to distance in $G$.
%$$E(H_{k,d}(G)) = \{ \{v_{1},v_2,\ldots, v_{t} \}: d_G(v_{j},v_{\ell}) \leq d \text{ for all } 1 \leq j \leq \ell \leq t \leq k,~\}.$$ 
This leads to the following theorem.

\begin{theorem}
Let $k \geq 2$ be a positive integer, and let $G$ be a graph. $G$ is $(k,d)$-special if and only if $H_{k,d}(G)$ has a perfect matching of cardinality $\left\lceil \frac{n}{k}\right\rceil$ such that at most one edge is not of order $k$.
\end{theorem}

\begin{proof}
Suppose $G$ is $(k,d)$-special and let $K=\{S_1,S_2,\ldots S_{\left\lceil \frac{n}{k} \right\rceil}\}$ be a witness of $G$ being $(k,d)$-special with $S_1$ being the set of order less than $k$ if necessary. It follows that $K$ is a set of disjoint edges in $H_{k,d}(G)$ whose union is the entire vertex set. Similarly, if $K$ is a perfect matching, then these sets have the desired cardinalities and the desired distances, so $G$ is $(k,d)$-special.
\end{proof}

Finally, we can bound $b_k$ of a graph using minimum edge covers in $H_{k,d}(G)$. An \emph{edge cover} in a hypergraph $H$ is a set of edges whose union is $V(H)$. Let $\beta'(H_{k,d}(G))$ be the minimum cardinality of an edge cover of $H_{k,d}(G)$.
    
\begin{theorem}\label{thm:edge_cover}
For positive integers $k,d$ and a graph $G$, we have that $$b_k(G) \leq \beta'(H_{k,d}(G)) + d.$$
Equivalently, if there exists a partition of $V(G)$ into $T$ sets each of order at most $k$ and diameter at most $d$ then$$b_k(G) \leq T + d.$$
\end{theorem}

Theorem~\ref{thm:edge_cover} will have application to bounds on $b_k$ for Cartesian grids; see Theorem~\ref{glast}.

\begin{proof}[Proof of Theorem~\ref{thm:edge_cover}]
Let $S_1,S_2,\dots,S_{\beta'(H_{k,d}(G))}$ be a minimum edge cover of \\ $H_{k,d}(G)$ and consider the $k$-liminal burning game on $G$. In each $S_i$, delete any vertices that are contained in some $S_j$ with $j > i$. Since this is a minimum cover, no $S_i$ is contained in the union of the other $S_j$, so after this process, all the $S_i$ are non-empty. We claim that in round $i$ for $1 \leq i \leq \beta'(H_{k,d}(G))$, the Arsonist can burn a vertex in some $S_i$ that does not yet contain a source. In round $i$, there have been $i-1$ sources placed so far, and the Saboteur has revealed $ki$ vertices. Since each $|S_i| \leq k,$ at most $k(i-1)$ vertices are in some $S_i$ that contains a source, so there are at least $k$ vertices in an $S_i$ that do not contain a source. After $\beta'(H_{k,d}(G))$ rounds, every $S_i$ contains a source, and after $d$ additional rounds, burning spreads to cover the remaining vertices since $\mathrm{diam}(S_i) \leq d.$

To show the equivalence with the second statement, observe that such a partition is a minimum edge cover for $H_{k,d}(G)$. Similarly, after deleting the overlap from the hyperedges, the $S_i$'s form a valid partition of $V(G)$.
\end{proof}

Theorem~\ref{thm:edge_cover} is a more general version of the bound given in Theorem~\ref{thm:kd_special_upperbound}, since being $(k,d)$-special implies the existence of a perfect matching of order $\left\lceil \frac{n}{k}\right\rceil$, which is also a vertex cover. The maximum matching problem is NP-complete~\cite{karp21}, and finding a minimum edge cover is an example of the set cover problem, which is also NP-complete~\cite{karp21}. The complexity of deciding if a graph is $(k,d)$-special remains an open problem.

\section{Cartesian products}
The \emph{Cartesian product} of graphs $G$ and $H$, denoted $G \square H$, has vertex set $V(G) \times V(H)$, and two distinct vertices $(u, v)$ and $(x, y)$ are adjacent in $G \square H$ if and only if $u = x$ and $(v, y) \in E(H)$, or $(u, x) \in E(G)$ and $v = y$. The notation $(V,h)$ will denote the set of vertices belonging to $V(G \Box H)$ that correspond to the set of vertices $V \subseteq V(G)$ and the vertex $h \in V(H)$. We define $(g,W)$ analogously. 

\subsection{Grids} \label{sec:grids}
The square \emph{Cartesian grid} graph $P_n \square P_n$ is denoted $G_{n}$. For square Cartesian grids, an asymptotically tight bound is known for the burning number \cite{bonf,MitschePralatRoshanbin}, and the cooling number is known up to an additive constant of two \cite{cooling}. 

\begin{figure}[htpb!]
\begin{tikzpicture}
  \begin{axis}[
      width=14cm, height=8cm,
      xlabel={$k$}, 
      xmin=0, xmax=600,
      ymin=0, ymax=200,
      grid=both,
    ]

    % Blue piecewise line
    \addplot[blue, thick, domain=0:50, samples=201] {200 - x/2 - 2};
    \addplot[blue, thick, domain=50:62, samples=13] {174};
    \addplot[blue, thick, domain=62:152, samples=344] {10000/x + 2*sqrt(x) - 3};
    \addplot[blue, thick, domain=152:600, samples=196] {10000/x + pow(3/2,1/3)*pow(100,2/3)};

    % Red piecewise line
    \addplot[red, thick, domain=0:223, samples=447] {200 - 3*x/4 - 7 - (2*x+3)*ln((2*100 + 2*x + 3)/(3*x + 7))};
    \addplot[red, thick, domain=223:405, samples=206] {10002/(x+2)};
    \addplot[red, thick, domain=405:600, samples=196] {pow(3/2,1/3)*pow(100,2/3) - 1};
    %\addlegendentry{Red}

    % Green data points
    \addplot[green, mark=*, only marks, mark size=1pt] coordinates {
      (1,189) (11,170) (21,161) (31,152) (41,146) (51,138)
      (61,128) (71,123) (81,116) (91,110) (101,104) (111,97)
      (121,92) (131,85) (141,80) (151,80) (161,73) (171,73)
      (181,68) (191,68) (201,67) (211,65) (221,65) (231,65)
      (241,60) (251,60) (261,54) (271,54) (281,54) (291,54)
      (301,54) (311,54) (321,54) (331,54) (341,53) (351,51)
      (361,51) %(371,51) (381,51) (391,51) (401,49) %(411,49)
      %(421,49) (431,49) (441,49) (451,49) (461,49) (471,49)
      %(481,49) (491,49) (501,49) (511,49) (521,49) (531,49)
      %(541,49) (551,49) (561,49) (571,49) (581,49) (591,49)
      %(601,49)
    };
    
    % Vertical dashed lines at k=50,200,405
    \addplot[dashed, gray, forget plot] coordinates {(50,0) (50,200)};
    \addplot[dashed, gray, forget plot] coordinates {(62,0) (62,200)};
    \addplot[dashed, gray, forget plot] coordinates {(152,0) (152,200)};
    \addplot[dashed, gray, forget plot] coordinates {(223,0) (223,200)};
    \addplot[dashed, gray, forget plot] coordinates {(405,0) (405,200)};

  \end{axis}
\end{tikzpicture}

\caption{The upper and lower bounds of $b_k(G_{100})$, along with a heuristic computational solution. The vertical dashed lines represent the transition from one range of $k$ to the next.}
\label{fig:grid_upper_lower100}
\end{figure}
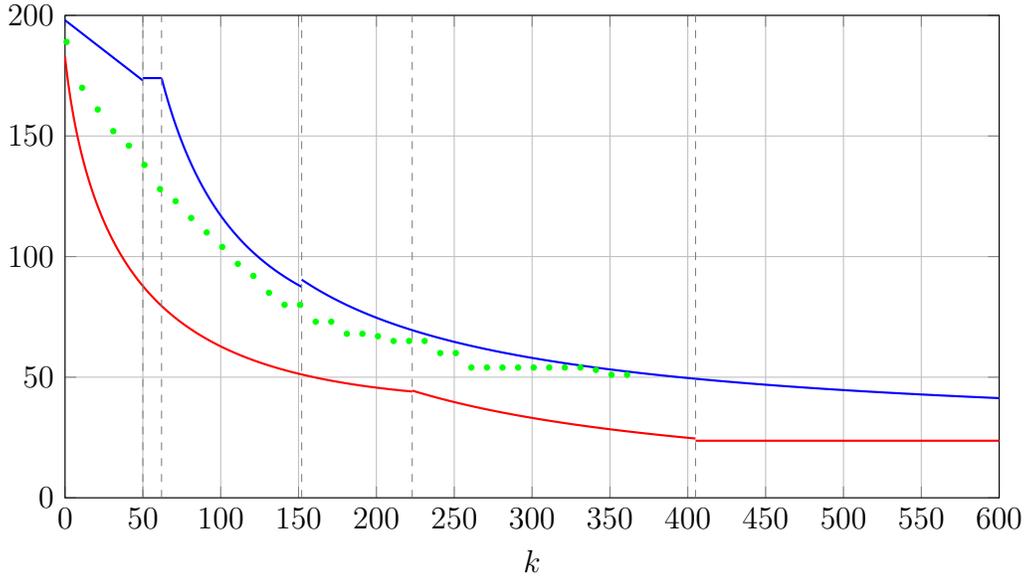

Here, we find $b_k(G_n)$ for all other ranges of $k$ up to constant factors, summarizing our results in the following table in seven cases for distinct ranges of $k,$ where $1\le k \le n^2$. 
In each case, we achieve the value of $b_k$ within a constant ratio, and provide the maximum ratio between the bounds for each range of $k$ as $n$ tends to infinity. We show these bounds for $n=100$ in Figure~\ref{fig:grid_upper_lower100}.

For ease of formatting in the table, we let $f(n) = 2n - 2\left\lfloor \log_2(n+3)\right\rfloor$, $g(k,n)=2n - \frac{3k}{4} -7 - (2k+3)\ln\left( \frac{2n+2k+3}{3k+7}\right)$, and $h(k,n)=\frac{n^2}{k} + (1+o(1))\left(\frac32\right)^{1/3}n^{2/3}$. 
We also let $n\geq 222$ so that $2n+23 < \left\lfloor \left(\frac{3}{16}\right)^{2/3} n^{4/3} \right\rfloor$ and the fifth interval is non-empty. 

\begin{table}[htpb!]
\centering
\begin{tabular}{c|c|c|c}
     Range of $k$& Lower & Upper & Max.\ Ratio\\%$O(b_k(G_n))$ \\
\hline
     $k=1$ & $f(n)$ & $f(n) +2$ & $1$\\%$f(n)+2$\\
     $2 \leq k < \frac{n-3}{2}$& $g(k,n)$ & $2n-\frac{k}{2} - 2$ &  $1.878$\\%$2n-\Theta(k)$\\
     $\frac{n-3}{2}\leq k \leq \frac{4n}{7}$& $g(k,n)$ & $\frac{7n}{4} -\frac{5}{4}$ &  $1.991$\\%$2n-\Theta(k)$\\
     $\frac{4n}{7} < k < 2n+23$ & $g(k,n)$ & $\frac{n^2}{k-2\sqrt{k}+1} + 2\sqrt{k} - 3$ & $2.146$\\ % $\Theta(n)$\\
     $2n+23 \leq k < \left(\frac{3}{16}\right)^{2/3} n^{4/3}$ & 
     $\frac{n^2+2}{k+2}$
     &$\frac{n^2}{k-2\sqrt{k}+1} + 2\sqrt{k} - 3$ 
     & $1$\\%$\Theta(\frac{n^2}{k})$\\
     $\left(\frac{3}{16}\right)^{2/3} n^{4/3} \leq k < \left(\frac{2}{3}\right)^{1/3} n^{4/3}$ & 
     $\frac{n^2+2}{k+2}$
     & $h(k,n)$
     & $2$\\%$\Theta(n^{2/3})$\\
     $k 
     \geq 
     \left(\frac23\right)^{1/3} 
     n^{4/3}$ 
     & $\left(\frac32\right)^{1/3}n^{2/3}-1$
     & $h(k,n)$
     & $2$\\% $\Theta(n^{2/3})$ \\
\end{tabular}
\caption{Values of $b_k(G_n)$ for various ranges of $k=k(n)$. Lower and upper bounds are provided along with the maximum ratio of the upper and lower bounds as $n$ tends to infinity.}
%asymptotics, up to highest order terms, of $b_k(G_n)$.}
\end{table}

For $0 \leq d \leq 2n-2$, let \emph{diagonal $d$} of $G_n$ be the vertices at distance exactly $d$ from $(0,0)$.
A vertex is \textit{lit} if it is burned or revealed, and \textit{unlit} otherwise.
We impose the \emph{graded-lexicographic} ordering on $V(G_n)$, where  $(a,b) \succ (c,d)$ if $a + b > c + d$ or $a + b = c + d$ and $a > c$. 
Define the \emph{basic strategy} for the Saboteur as the strategy where the Saboteur always reveals the $k$ smallest vertices that are unburned. 
The basic strategy works well for the Saboteur, and we use it to obtain the improved lower bounds. 
When $k$ is not too large, a sensible response for the Arsonist when playing against this basic strategy is to play the smallest revealed vertex.

We also recall here that the length of the game that occurs with a fixed Saboteur strategy is never longer 
than when the Saboteur plays optimally. Thus, the length of the game with a fixed Saboteur strategy and the Arsonist playing optimally is a lower bound on the liminal burning number for the graph. 

\begin{theorem}\label{gridLB}
For positive integers $k$ and $n,$ we have the following.
\[b_k(G_n) \geq 
\begin{cases}
        2n - 2\left\lfloor \log_2(n+3)\right\rfloor
            & \text{for $k=1$,}\\
		
		2n - \frac{3k}{4} -7 - (2k+3)\ln\left( \frac{2n+2k+3}{3k+7}\right) 
            & \text{for $2 \leq k < 2n+23$,}\\ 
        
        \frac{n^2+2}{k+2} 
            & \text{for $2n+23 \leq k < \left(\frac23\right)^{1/3} n^{4/3}$,}\\
        
        \left(\frac32\right)^{1/3}n^{2/3}-1
            & \text{for $k \geq \left(\frac23\right)^{1/3} n^{4/3}$.}\\
		 \end{cases}
\]
\end{theorem}
\begin{proof}
    Liminal burning when $k=1$ is equivalent to the cooling process \cite{cooling} and yields the given result. 
    In the case that $k \geq \left(\frac23\right)^{1/3} n^{4/3}$, it is known that $b(G_n) \geq \left(\frac32\right)^{1/3}n^{2/3}-1$ \cite{MitschePralatRoshanbin}, so the result follows since $b(G_n) \leq b_k(G_n) $. 

    For the remaining cases, we will study this game over four phases. 
    During each of these phases, the Saboteur is forced to use the basic strategy. 
    
    \emph{\textbf{Phase 1:}}
    For a recursion, at the start of round $i$, we assume that there is at least one vertex $v$ on diagonal $d$ such that $v$ has one or two neighbors which are unlit, and such that any vertex $u$ with $u \prec v$ is lit. 
    This is true in the base case of $i=2$, since the Saboteur reveals the first $k$ vertices in the ordering in the first round, and any of these vertices that the Arsonist chooses to burn is adjacent to at most two unrevealed vertices.

    For the general case, we play while $d \leq \left\lfloor \frac{k-1}{2} \right\rfloor$. 
    The vertices adjacent to $v$ that were unlit are burned by propagation, and now it is the Saboteur's move. 
    Between vertex $v$ and the pair of most recently burned vertices, there are at most $d-1$ unrevealed vertices in the ordering, which the Saboteur now reveals.  
    The next $d+2$ in the ordering that are not yet revealed are exactly those vertices that are adjacent to a vertex that is already lit. 
    As $k \geq (d+2)+(d-1)$, the Saboteur reveals these vertices, and so now only unburned vertices are adjacent to the vertices that are unlit. 
    The Arsonist burns a revealed vertex, which is now the only possible burned vertex next to a vertex that is unlit. Since each revealed vertex that the Arsonist may have chosen is adjacent to at most two vertices that are unlit, the general recursive step holds.

    Note that $k$ vertices were revealed on the first round and $k+2$ vertices were lit on each subsequent round during the recursion, yielding $(k+2)i -2$ vertices lit at the end of round $i$. 
    Since we revealed at least one vertex in diagonal $\left\lfloor \frac{k-1}{2} \right\rfloor+2$, all vertices up to and including the vertices in diagonal $\left\lfloor \frac{k-1}{2} \right\rfloor+1$ were lit, meaning $(k+2)i -2 \geq \left(\genfrac{}{}{0pt}{}{\left\lfloor \frac{k-1}{2} \right\rfloor+2}{2}\right)$ for some round $i$ during this process. Thus, $i\geq \frac{k}{8}$, and so at least $\left\lceil \frac{k}{8} \right\rceil$ rounds have elapsed and we have not lit a vertex in diagonal $\left\lfloor \frac{k-1}{2} \right\rfloor+3 = \left\lfloor \frac{k+5}{2} \right\rfloor$ or later. 

    Note in the case that $\left\lfloor \frac{k+5}{2} \right\rfloor\geq n$, we can continue to burn or reveal $k+2$ vertices during each round except the first. 
    This means that in round $\left\lceil \frac{n^2+2}{k+2} \right\rceil$, we reveal or burn a vertex that was unlit. 

\emph{\textbf{Phase 2:}}
    Define $T_1 = \left\lceil \frac{k}{8} \right\rceil$ and $d_1 = \left\lfloor \frac{k+5}{2} \right\rfloor$.
    We proceed to phase 2, starting in round $T_1 +1$ with all vertices that are lit strictly before diagonal $d_1$. 
    We begin by burning any unburned vertex strictly before diagonal $d_1$. 

    During each round of this phase, we burn all revealed vertices during the Arsonist's move, so exactly $k$ additional vertices will be burned each round in addition to the vertices burned through propagation. 

    Let $V_i$ denote the number of burned vertices at the end of round $T_1+i$, where $V_0 = \binom{d_1}{2}$. 
    During each round, propagation ensures that each vertex that has a burned neighbor is itself burned. 
    This yields $p_i$ new burned vertices that are burned during round $i$ if $\binom{p_i-1}{2}+1 \leq V_{i-1} \leq \binom{p_i}{2}$, which implies that $p_i \leq \sqrt{2V_{i-1}-\frac{7}{4}}+\frac{3}{2}$.
    After propagation, the Saboteur reveals $k$ vertices, which are burned by the Arsonist. 
    This implies that the total number of added burned vertices in round $i$ was $V_i-V_{i-1} = p_i +k \leq \sqrt{2V_{i-1}-\frac{7}{4}}+\frac{3}{2} +k$. 
    Let $f:\mathbb{N} \rightarrow \mathbb{R}$ by $f(0)=V_0$ and $f'(x)=\sqrt{2f(x)-\frac{7}{4}}+\frac{3}{2} +k$.
    We will show that $f(i)\geq V_i$ for all $i$. 
    
    For an induction, assume that $f(i)\geq V_i$ and $f'(i)\geq V_{i+1}-V_{i}$. 
    The base case where $i=0$ is immediate. 
    Since $f''(x)>0$ and $f'(i)\geq V_{i+1}-V_i$, we see that $f(i+1) \geq V_i + (V_{i+1}-V_i) = V_{i+1}$, showing the first inequality. 
    By definition, 
    \begin{equation*}
    f'(i+1) = \sqrt{2f(i+1)-\frac{7}{4}}+\frac{3}{2} +k 
    \geq \sqrt{2V_{i+1}-\frac{7}{4}}+\frac{3}{2} +k 
    \geq p_{i+2} + k = V_{i+2}-V_{i+1},
    \end{equation*}
    showing the second inequality and completing the induction. 
    
    Thus $f(i)$ gives an upper bound on the number of burned vertices at time $T_1+i$. 
    We must find a good estimate of the number of rounds that can be played while ensuring that no vertex in diagonal $n$ or later is burned.
    If $f(T)\leq \binom{n}{2}$, then  $V_T \leq f(T) \leq \binom{n}{2}$, so all burned vertices will be in a diagonal less than $n$. 
    Hence, $T_1+\lfloor{T}\rfloor$ is our estimate, and we proceed by calculating $T$. 
    
    Define $g(x) = \sqrt{2f(x)-\frac{7}{4}}$, noting that $f(x) = \frac{1}{2}\left(g(x)^2 + \frac{7}{4}\right)$ as a result. 
    Taking the derivative of both sides, we see that $f'(x) = g(x)g'(x)$. However, by definition, we also have $f'(x) = \sqrt{2f(x)-\frac{7}{4}}+\frac{3}{2} +k = g(x) + \frac{3}{2}$, implying that $g(x)g'(x) = g(x) +\frac{3}{2} +k$, from which $g'(x) = 1 + \left(\frac{3}{2}+k\right)g(x)^{-1}$. Separating variables yields $dt = g\left(g + \frac{3}{2}+k\right)^{-1}dg$.
    Note that $g(0)=\sqrt{d_1^2-d_1 -\frac{7}{4}}$ and $g(T)=\sqrt{n^2-n -\frac{7}{4}}$. 
    We integrate $dt$ from $x=0$ to $x=T$ to find 
    \begin{align*} 
    T &= \int_{g(0)}^{g(T)} 1+\frac{k+\frac{3}{2}}{g} dg =  \left. g - \left(k + \frac{3}{2}\right) \ln \left(g+\frac{3}{2} +k\right)\right|^{g=\sqrt{n^2-n -\frac{7}{4}}}_{g=\sqrt{d_1^2-d_1 -\frac{7}{4}}}\\
    &= \sqrt{n^2-n-\frac{7}{4}} - \sqrt{d_1^2-d_1-\frac{7}{4}} - \left(k+\frac{3}{2}\right)\ln \left({\frac{\sqrt{n^2 - n-\frac{7}{4}}+\frac{3}{2}+k}{\sqrt{d_1^2 - d_1-\frac{7}{4}}+\frac{3}{2}+k}}\right).
    \end{align*}
    The mean value theorem shows that $\sqrt{n^2-n-\frac{7}{4}} - \sqrt{d_1^2-d_1-\frac{7}{4}} \geq n-d_1$, and showing that $\frac{\sqrt{x^2 - x-\frac{7}{4}}+\frac{3}{2}+k}{x+\frac{3}{2}+k}$ is monotonically increasing is straightforward, which implies that     
    ${\frac{\sqrt{n^2 - n-\frac{7}{4}}+\frac{3}{2}+k}{\sqrt{d_1^2 - d_1-\frac{7}{4}}+\frac{3}{2}+k}} \geq {\frac{n+\frac{3}{2}+k}{d_1+\frac{3}{2}+k}}$. 
    Therefore, after round 
    $$T_1+\left\lfloor n-d_1-\left(k+\frac{3}{2}\right) \ln \left({\frac{n+\frac{3}{2}+k}{d_1+\frac{3}{2}+k}}\right) \right\rfloor,$$ no vertex in diagonal $n$ or later is lit. 
    
\emph{\textbf{Phase 3:}}
We define $T_2 = \left\lfloor n-d_1-\left(k+\frac32\right) \ln \left({\frac{n+\frac32+k}{d_1+\frac32+k}}\right) \right\rfloor$, so that $T_1+T_2$ rounds have elapsed at the start of Phase 3. 
As in Phase 2, we burn all revealed vertices during the Arsonist's move. 

Let $S_i$ be the set of vertices that were burned in round $T_1+i$, for $1 \leq i \leq T_2$. 
We note that these sets satisfy the following three properties:
\begin{enumerate}
\item the distance from a vertex in $S_i$ to diagonal $d_1$ was at least $i$;  
\item $d(S_i,S_j) \geq j-i+1$, where $i<j$; and 
\item the distance from a vertex in $S_i$ to diagonal $n$ was at least $T_2+1-i$. 
\end{enumerate}
These conditions are required for $S_i$ to be unlit at the start of round $T_1+i$, and for no vertex in diagonal $n$ or more to be lit by round $T_1+T_2$. 

Define $S'_i = \{(n-1-x,n-1-y) : (x,y)\in S_{T_2+1-i}\}$ for $1 \leq i \leq T_2$. 
Inherited from the properties of $S_i$ are the the following properties for $S'_i$:

\begin{enumerate}
\item the distance from a vertex in $S'_i$ to a vertex of distance $d_1$ from $(n-1,n-1)$  was at least $i$;
%\item the distance from a vertex in $S'_i$ to diagonal $2n-2-d_1$ was at least $i$;
\item $d(S'_i,S'_j) \geq j-i+1$, where $i<j$; and 
%\item the distance from a vertex in $S'_i$ to diagonal $n-2$ was at least $T_2+1-i$.
\item the distance from a vertex in $S'_i$ to a vertex of distance $n$ from $(n-1,n-1)$ was at least $T_2+1-i$. 
\end{enumerate}
The Saboteur can then reveal the $k$ vertices in $S'_i$ at time $T_1+T_2+i$, which the Arsonist immediately burns. 
At the end of round $T_1+T_2+T_2$, this leaves the vertices of distance at most $d_1$ from $(n-1,n-1)$ both unlit.   
    
\emph{\textbf{Phase 4:}}
To begin, we burn any unburned vertex of distance more than $d_1$ from $(n-1,n-1)$. 
    During the first propagation of this phase, the vertices of distance exactly $d_1$ from $(n-1,n-1)$ are burned. 
    This leaves exactly $\binom{d_1}{2}$ unlit vertices. 
    In round $T_1+2T_2+1$, we have that $k$ of these vertices will be revealed, one of which will be burned by the Arsonist. 
    Imitating the first phase, $k+2$ vertices will be lit on each subsequent round. 
    If $T_4$ rounds are performed in this phase, then $(k+2)T_4-2$ vertices are lit. Since $(k+2)(T_1-1)-2 < \binom{d_1}{2}$, at least $T_1$ rounds are required in this phase to allow all vertices to be lit. 
    
    Thus at least $2(T_1+T_2)$ rounds must be played for all of $G_n$ to be burned. Considering the parity of $k$ in $d_1$ and simplifying, we find that $2T_1 \geq \frac{k}{4}$ and $$2T_2 \geq 2n-k-7-(2k+3)\ln \left({\frac{2n+2k+3}{3k+7}}\right),$$ and the proof follows. 
\end{proof}

We also have upper bounds on $b_k(G_n)$.
\begin{theorem}\label{grub}
For positive integers $k$ and $n$, we have the following.
\[b_k(G_n) \leq 
\begin{cases}
        2n - 2\left\lfloor \log_2(n+3)\right\rfloor+2
            & \text{for $k=1$,}\\
		
		2n - \frac{k}{2} -2 
            & \text{for $2 \leq k \leq \frac{n-3}{2}$,}\\
            
		\frac{7n}{4} -\frac{5}{4} 
            & \text{for $ \frac{n-3}{2} < k \leq \frac{4n}{7}$,}\\
        
        \frac{n^2}{k-2\sqrt{k} + 1} + 2\sqrt{k}-3 
            & \text{for $\frac{4n}{7} < k < \left(\frac{3}{16}\right)^{2/3} n^{4/3}$,}\\
        
        \frac{n^2}{k} + (1+o(1))\left(\frac32\right)^{1/3}n^{2/3} 
            & \text{for $k \geq \left(\frac{3}{16}\right)^{2/3} n^{4/3}$.}\\
		 \end{cases}
   \] 
\end{theorem}

The result for $k=1$ follows from \cite{cooling}.
    For $k \geq \left(\frac{3}{16}\right)^{2/3} n^{4/3}$, it is known that $b(G_n) \leq (1+o(1))\left(\frac32\right)^{1/3}n^{2/3}$ \cite{MitschePralatRoshanbin}. The Arsonist can wait for the first $\frac{n^2}{k}$ rounds while all vertices are revealed, then burn the graph in an additional $(1+o(1))\left(\frac32\right)^{1/3}n^{2/3}$ rounds. 
    The case with $ \frac{n-3}{2} < k \leq \frac{4n}{7}$ follows since $b_k(G_n) \leq \frac{7n}{4} -\frac{5}{4}$ when $k=\frac{n-3}{2}$ and since $b_x(G_n)\geq b_y(G_n)$ for $x\leq y$. The remaining two cases of Theorem~\ref{grub} are separated into theorems as follows. We first consider $k \leq \frac{n-3}{2}$.

\begin{theorem}
For positive integers $k$ and $n$ with $2\leq k\leq \frac{n-3}{2}$, we have that
\[
    b_k(G_n) \leq 2n-\frac{k}{2} - 2 .
\]
\end{theorem}
\begin{proof}
We play a game of liminal burning on $G_n$, but we may force the Arsonist to play a certain strategy, or disallow the Arsonist from playing their move, or allow the Saboteur to play less than $k$ vertices. 
Each of these actions can only increase the number of rounds. 
We show that all vertices will be burned after round $\left\lfloor 2n - \frac{k}{2} - \frac52 \right\rfloor$, independent of the strategy of the Saboteur. 

    However, the Saboteur and Arsonist play their first moves, we may rotate the grid $G_n$ so that the Arsonist's first burned vertex is closer to the vertex $(0,0)$ than to any other corner vertex. 
    The only vertices that can be of distance more than $2n - \frac{k}{2} - \frac52 - 1$ from this first vertex are the vertices of distance at most $\frac{k}{2} + \frac52 - 1$ from the vertex $(n-1,n-1)$. 
    As a result, after $2n - \frac{k}{2} - \frac52 - 1$ propagation rounds, any unburned vertex must be in this set of vertices. 
    
    On each round $t$ with $1 \leq t \leq \left\lfloor \frac{k}{2} \right\rfloor$, we force the Arsonist to play on a revealed vertex with a neighbor that is unlit, where possible. 
    During the following round, since there must be a burned vertex next to an unlit vertex, at least one vertex will be burned during propagation. 
    A further $k$ vertices will be revealed, so at least $k+1$ vertices that were unlit are lit at the end of this next round. 
    Thus, at the end of round $T_1 = \left\lfloor\frac{k+2}{2} \right\rfloor$, at least $(k+1)\left\lfloor\frac{k+2}{2}\right\rfloor-1 \geq \binom{k+1}{2}+1$ vertices are lit.
     
At least one vertex of distance $k+1$ or more from $(0,0)$, say $(x,y)$, must have been burned during or before the end of round $T_1$.
If $x,y\leq \left\lfloor \frac{n}{2} \right\rfloor$, then $(x,y)$ has distance at most $2n-2-(k+1)$ to any $(x',y')$ with $x',y' > \left\lfloor \frac{n}{2} \right\rfloor$. As such, all vertices would be burned by the end of round $T_1+2n-2-(k+1) \leq 2n - 2 - \frac{k}{2}$. 
If either $x$ or $y$ is larger than $\left\lfloor \frac{n}{2} \right\rfloor$, then $(x,y)$ has distance at most $\frac{3}{2}n - \frac32$ to any $(x',y')$ with $x',y' > \left\lfloor \frac{n}{2} \right\rfloor$. As such, all vertices would be burned by the end of round $T_1+\frac{3}{2}n - \frac32 \leq \frac{3}{2}n + \frac{k}{2} - \frac12 \leq 2n - \frac{k}{2} - 2$ since $k\leq \frac{n-3}{2}$.
\end{proof}

We now consider the final range of $k$ by using Theorem~\ref{thm:edge_cover}.

\begin{theorem}\label{glast}
For integers $k,n\geq 1$, 
    \[
    b_k(G_n) \leq \frac{n^2}{k-2\sqrt{k}+1} + 2\sqrt{k} - 3.
    \]
\end{theorem}
\begin{proof}
    Define $K=\left\lfloor \sqrt{k} \right\rfloor. $
    We decompose the vertices $V(G_n)$ into the subsets $$V_{i,j} = \{(K\cdot i + x, K\cdot j + y) : 0 \leq x,y \leq K\}$$ for $0 \leq i,j \leq \left\lceil \frac{n}{K} \right\rceil$. 
    While most of these subsets contain $K^2 \leq k$ vertices, some may contain fewer if $i$ or $j$ equal $\left\lceil \frac{n}{K} \right\rceil$. This partition has $\left\lceil \frac{n}{K} \right\rceil ^2$ sets, each with cardinality at most $k$ and diameter at most $2K - 2$. Hence, by Theorem~\ref{thm:edge_cover} we have that
    \[
    b_k(G_n) \leq \left(\left\lceil \frac{n}{K} \right\rceil \right)^2 + 2K - 2 \leq \frac{n^2}{k-2\sqrt{k}+1} + 2\sqrt{k} - 3.
    \]
The proof follows.
\end{proof}

Theorem~\ref{glast} implies that Theorem~\ref{thm:edge_cover} is asymptotically tight for $2n+23 \leq k < \left(\frac{3}{16}\right)^{2/3} n^{4/3}$ since the ratio between our upper and lower bounds tends to $1$.

\subsection{General Cartesian Products}
We now turn to liminal burning in more general Cartesian products. In which graphs is $b_k$ equal to the cooling number for $k>1$? One approach is to have each vertex of the original graph behave like a $k$-clique so that Saboteur reveals only one ``vertex'' per round, but this approach does not always yield cooling behaviour. However, if a large enough clique replaces each vertex, the two parameters become equal. 

\begin{theorem}\label{thm:bk_eqls_cooling}
Let $G$ be a graph, and let $j,k$ be positive integers. If $j \geq \mathrm{diam}(G)+k$, then $b_k(G \Box K_j) = \mathrm{CL}(G \Box K_j)$.  
\end{theorem}

\begin{proof}
Consider a cooling sequence of $G \Box K_j$ given by $(c_1,c_2,\ldots,c_m)$ and partition $V(G \Box K_j)$ into sets $(g,V(K_j))$ for each $g \in V(G)$. Each $c_i$ has a unique corresponding vertex set $V_i =(c_i,V(K_j))$, as each $V_i$ becomes fully burned after the selection of $c_i$ by the Arsonist. We propose the following strategy for the Saboteur: In round $i$, the Saboteur reveals $k$ arbitrarily chosen vertices within $V_i$. The Arsonist must then select a vertex within $V_i$ to burn, as all vertices in $V_1,V_2,\ldots,V_{i-1}$ are burned by this round. Each previously chosen source by the Arsonist produces at most 1 burned vertex in $V_i$ by round $i$ (its ``image'' in $V_i$), so as $|V_i| = j \geq \mathrm{diam}(G)+k$, the Saboteur will always have at least $k$ vertices to reveal each round. Thus, the Arsonist must burn a vertex in $V_i$ in round $i$ for $1 \leq i \leq m$, ensuring the game lasts at least $\mathrm{CL}(G \Box K_j)$ rounds. 
\end{proof}

This idea of revealing $k$ vertices in $G \Box K_j$ in a way that mimics revealing only one vertex in $G$ can be carried through to the more general setting of $G \Box H$. If $H$ is large enough compared to $G$, then the Saboteur can force the Arsonist to cool $G$. If $H$ is even larger, then the faux cooling performed on $G$ becomes optimal.

\begin{theorem}\label{thm:GxH_lowerbounds}
Let $G, H$ be graphs, and let $k$ be a positive integer. The following statements are true.
    
\begin{enumerate}
\item If $|V(H)| \geq k$, then $b_k(G \Box H) \geq \mathrm{CL}(G)$.
\item If $|V(H)| \geq \mathrm{diam}(G)+k$, then $b_k(G \Box H) \geq \mathrm{diam}(G)+1$.
\end{enumerate}
    
\end{theorem}

\begin{proof}
Partition $V(G\Box H)$ into sets $(g,V(H))$ for each $g \in V(G)$. For (1), let $c_1,c_2,\ldots,c_m$ be a cooling sequence for $G$, and note that each $c_i$ has a unique corresponding vertex set $V_i = (c_i,V(H))$. The strategy for the Saboteur is as follows: In round $i$, reveal $k$ vertices chosen arbitrarily in $V_i$. Each $V_i$ contains at least $k$ vertices, and no previously burned vertex produces a burned vertex in $V_i$ during this round because $d(V_i,V_j) = d(c_i,c_j) \geq |i-j|+1$. Thus, during round $ 1 \leq i \leq m$, all of $V_i$ is unlit, ensuring the game lasts at least $\mathrm{CL}(G)$ rounds.

For (2), let $p_1,p_2,\ldots,p_d$ be a diameter length path in $G$, and observe that each $p_i$ has a unique corresponding vertex set $V_i = (p_i,V(H)).$ The strategy for the Saboteur is as follows: In round $i$, the Saboteur will reveal $k$ vertices in $V_i$. This is always possible, as in round $i$ the vertex set $V_i$ will contain at most $i-1$ burned vertices from the Arsonist's previously chosen sources, and as $|V_i| = \mathrm{diam}(G)+k \geq 1+k-1$ for $1 \leq i \leq \mathrm{diam}(G)+1$, the game must last at least $\mathrm{diam}(G)+1$ rounds.
\end{proof}

We extend this cooling-type argument to multiple Cartesian products of a graph $G$, stated in terms of the length of a maximum cooling sequence in $G$ as opposed to $\mathrm{CL}(G)$. This is because the sets that are revealed as part of the $k$-liminal burning game are $k$-subsets of the cooling sequence of $G$ and therefore cannot reflect the behavior of $\mathrm{CL}(G)$ during the last round of cooling $G$. For a positive integer $n$ and a graph $G$, let $G^{\Box n}$ denote the Cartesian product of $n$ copies of $G$.

\begin{theorem}\label{thm:bk_Gn}
Let $k,n$ be positive integers and let $G$ be a graph. If $m$ is the length of an optimal cooling sequence of $G$, then $$b_k\left(G^{\Box n}\right) \geq \left\lfloor \frac{m}{k} \right\rfloor + (n-1)(m-1).$$
\end{theorem}

\begin{proof}
Let $c_1,c_2,\ldots,c_m$ be a cooling sequence for $G$ and let $$V_i = \{c_{(i-1)k+1},\ldots,c_{(i-1)k+k}\}.$$ for $1 \leq i \leq \left\lfloor \frac{m}{k} \right\rfloor$. The strategy for the Saboteur is as follows: In round $i$, reveal the set of vertices $S_i$, where $S_i = ((V_i)_1,(c_1)_2,(c_1)_3,\ldots,(c_1)_m)$ for $1 \leq i \leq \left\lfloor \frac{m}{k} \right\rfloor,$ and for $i = \left\lfloor \frac{m}{k} \right\rfloor+(m-1)(q-2)+1$ where $2 \leq q \leq n$ and $1 \leq j \leq m-1,$ let
\begin{equation*}
    S_i =   \left(\left(V_{\left\lfloor \frac{m}{k}\right\rfloor}\right)_1,\left(c_{m}\right)_2,\left(c_m\right)_3,\ldots,\left(c_m\right)_{q-1},\left(c_{j}\right)_q,\left(c_1\right)_{q+1},\ldots,\left(c_1\right)_m\right). 
\end{equation*}
We note that $d(S_i, S_j) \geq |i-j|+1$ for $1\leq  i \neq j \leq \left\lfloor \frac{m}{k} \right\rfloor + (n-1)(m-1)$, showing that no matter the vertex choices of the Arsonist, the Saboteur can reveal these vertices each round as planned. Therefore, the game must last at least as long as the Saboteur's strategy takes and $b_k\left(G^{\Box n}\right) \geq \left\lfloor \frac{m}{k} \right\rfloor + (n-1)(m-1).$
\end{proof}

The following corollary is immediate.

\begin{corollary}
Let $n$ be a positive integer and $G$ a graph. If $m$ is the length of an optimal cooling sequence of $G$, then $\mathrm{CL}\left(G^{\Box n}\right) \geq n(m-1)+1.$
\end{corollary}

\section{Complexity}

We next investigate the complexity of determining the length of the $k$-liminal burning game. We recommend~\cite{AroraBorak} for a background on computational complexity. To proceed, we define the following decision problem for $k=f(n)$.
    
    \medskip
    
    \noindent PROBLEM: $k$-\textsc{Liminal-Burning}.\\
    INSTANCE: A graph $G$, and an integer $c$.\\
    QUESTION: Is $c \leq b_k(G)$?
    
    \medskip

We choose $c \leq b_k(G)$ rather than $c = b_k(G)$ or $c \geq b_k(G)$ as this is analogous to asking if the Saboteur ``wins the game.'' For a clear winner, we could add a fixed parameter $c$ and say the Saboteur wins if the game lasts at least $c$ rounds. 
    
As instances of $k$-\textsc{Liminal-Burning}, {\sc Burning} ($k=|V(G)|$) and {\sc Cooling} ($k=1$) are co-NP-complete \cite{Burning_Hard, Burning_APX} and NP-complete \cite{Cooling_Complexity} decision problems for graphs respectively. As a two-player game, liminal burning is more comparable to combinatorial games like Go~\cite{PSPACE_Go} and Generalized Geography~\cite{PSPACE_GG}. 

We show that $k$-\textsc{Liminal-Burning} is PSPACE-complete, first defining the PSPACE-complete problem {\sc $3$-QBF} \cite{PSPACE_QBF}. Recall that a $3$-CNF formula is made up of clauses $x \lor y \lor z$ which are all joined by conjunctions.

    \medskip

    \noindent PROBLEM: {\sc $3$-QBF}.\\
    INSTANCE: A quantified boolean formula $Q_1 x_1 Q_2 x_2\cdots Q_n x_n\phi(x_1,\dots,x_n)$ with $\phi$ a $3$-CNF formula.\\
    QUESTION: Is $Q_1 x_1Q_2 x_2 \cdots Q_n x_n\phi(x_1,\dots,x_n)$ true?\\
    
    \medskip
    
The majority of the section will be devoted to constructing a reduction from {\sc $3$-QBF} to $k$-\textsc{Liminal-Burning}.

\begin{theorem} \label{thm:pspace}
For all integers $k\geq 2$, $k$-{\sc Liminal-Burning} is PSPACE-complete.
\end{theorem}

We first consider a construction based on the one in~\cite{Cooling_Complexity}. Given a $3$-CNF formula $\phi$ with $n$ variables and $m$ clauses, we construct a graph $G_\phi$ as follows. First, create a path $a_0,a_1,\dots a_{n}$. For each variable $x_i$ with $1 \leq i \leq n,$ add vertices $x_i$ and $\lnot x_i$, add the edges $x_na_n$ and $\lnot x_na_n$, and for $1 \leq i \leq n-1$ add edges $x_i\lnot x_i$ $x_ia_i$, $\lnot x_ia_i$, $x_ia_{i+1}$, and $\lnot x_ia_{i+1}$. Name these $x_i$ \emph{variable vertices}. For each clause $C_j$, add a copy of $K_3$ whose vertices correspond to the literals in $C_j$ called a \emph{clause gadget}. If $x_i \in C_j$, then add a path with $n + j - i - 1$ internal vertices between $x_i$ and its corresponding vertex in the $C_j$ clause gadget. 

The idea of the construction is that if $a_0$ is burning before the game starts and $\phi$ is satisfiable, then given a satisfying assignment for $\phi$, cooling whichever of $x_i,\lnot x_i$ is false in round $i$ will cool the graph in at least $n + m$ rounds, since in round $n + j$  the $C_j$ clause gadget which has not yet been cooled. See Figure~\ref{fig:Gphi} for an example.

    \begin{figure}[htpb!]
            \centering 
            \begin{tikzpicture}
                \GraphInit[vstyle=Classic]
                \SetUpVertex[FillColor=white]
            
                \tikzset{VertexStyle/.append style={minimum size=24pt, inner sep=1pt}}
            
                \Vertex[x=0,y=0,LabelOut=false,L={$a_0$},]{A0}
                \Vertex[x=2,y=0,LabelOut=false,L={$a_1$},]{A1}
                \Vertex[x=4,y=0,LabelOut=false,L={$a_2$},]{A2}
                \Vertex[x=6,y=0,LabelOut=false,L={$a_3$},]{A3}

                \Vertex[x=3,y=1,LabelOut=false,L={$x_1$},]{X1}
                \Vertex[x=5,y=1,LabelOut=false,L={$x_2$},]{X2}
                \Vertex[x=7,y=1,LabelOut=false,L={$x_3$},]{X3}

                \Vertex[x=3,y=-1,LabelOut=false,L={$\lnot x_1$},]{NX1}
                \Vertex[x=5,y=-1,LabelOut=false,L={$\lnot x_2$},]{NX2}
                \Vertex[x=7,y=-1,LabelOut=false,L={$\lnot x_3$},]{NX3}

                % Clause 1

                \Vertex[x=1,y=3,LabelOut=false,L={$x_1$},]{C1X1}
                \Vertex[x=2,y=4,LabelOut=false,L={$x_2$},]{C1X2}
                \Vertex[x=3,y=3,LabelOut=false,L={$\lnot x_3$},]{C1NX3}
                
                % Clause 2

                \Vertex[x=5,y=3,LabelOut=false,L={$ \lnot x_1$},]{C2NX1}
                \Vertex[x=6,y=4,LabelOut=false,L={$\lnot x_2$},]{C2NX2}
                \Vertex[x=7,y=3,LabelOut=false,L={$x_3$},]{C2X3}

                \tikzset{VertexStyle/.append style={minimum size=12pt, inner sep=1pt}}

                % C1X1 path
                
                \Vertex[x=1,y=2,NoLabel=true,]{C1X1B}
                \Vertex[x=1.5,y=1.5,,NoLabel=true,]{C1X1A}

                % C1X2 path

                \Vertex[x=2.3,y=2,NoLabel=true,]{C1X2A}

                % C2NX1 path
                
                \Vertex[x=3.8,y=1,NoLabel=true,]{C2NX1A}
                \Vertex[x=4,y=2,,NoLabel=true,]{C2NX1B}
                \Vertex[x=4,y=3,,NoLabel=true,]{C2NX1C}
                
                % C2NX2 path

                \Vertex[x=6,y=2,NoLabel=true,]{C2NX2A}
                \Vertex[x=6,y=2.7,NoLabel=true,]{C2NX2B}

                % C2X3 path

                \Vertex[x=7,y=2,NoLabel=true,]{C2X3A}

                \tikzset{VertexStyle/.append style={draw=white, minimum size=24pt}}
                
                \Vertex[x=2,y=5,LabelOut=false,L={$C_1$}]{C1}
                \Vertex[x=6,y=5,LabelOut=false,L={$C_2$}]{C2}
                
                \Edge(A0)(A1)
                \Edge(A1)(A2)
                \Edge(A2)(A3)

                \Edge(X1)(NX1)
                \Edge(X2)(NX2)
                \Edge(X3)(NX3)

                \Edge(X1)(A1)
                \Edge(A1)(NX1)
                \Edge(X1)(A2)
                \Edge(A2)(NX1)
                
                \Edge(X2)(A2)
                \Edge(A2)(NX2)
                \Edge(X2)(A3)
                \Edge(A3)(NX2)
                
                \Edge(X3)(A3)
                \Edge(A3)(NX3)

                \Edge(C1X1)(C1X2)
                \Edge(C1X1)(C1NX3)
                \Edge(C1NX3)(C1X2)

                \Edge(C2NX1)(C2NX2)
                \Edge(C2NX1)(C2X3)
                \Edge(C2X3)(C2NX2)

                \Edge(X1)(C1X1A)
                \Edge(C1X1B)(C1X1A)
                \Edge(C1X1B)(C1X1)

                \Edge(X2)(C1X2A)
                \Edge(C1X2A)(C1X2)

                % C2NX1 path
                
                \Edge(NX1)(C2NX1A)
                \Edge(C2NX1A)(C2NX1B)
                \Edge(C2NX1B)(C2NX1C) 
                \Edge(C2NX1C)(C2NX1)
                
                % C2NX2 path

                \Edge(NX2)(C2NX2A)
                \Edge(C2NX2A)(C2NX2B)
                \Edge(C2NX2B)(C2NX2)

                % C2X3 path

                \Edge(X3)(C2X3A)
                \Edge(C2X3A)(C2X3)
                
                \tikzset{EdgeStyle/.style = {bend right=20}}

                \Edge(NX3)(C1NX3) % bend
                
            \end{tikzpicture}

            \caption{The graph $G_\phi$ when $\phi(x_1,x_2,x_3)=(x_1 \lor x_2 \lor \lnot x_3) \land (\lnot x_1 \lor \lnot x_2 \lor x_3)$.}
            \label{fig:Gphi}
        \end{figure}

This construction provides a mechanism for translating a satisfying assignment for a Boolean formula into a cooling sequence where the first vertex of the sequence is fixed.
To complete our reduction, we need to do three things. First, we simulate quantifiers on the variables of the formula, so the Saboteur assigns values to existentially quantified variables and the Arsonist assigns values to universally quantified variables. To do so, we create $k$ identical vertices in place of any variable vertices which were existentially quantified and use $\left\lceil \frac{k}{2}\right\rceil$ and $\left\lfloor \frac{k}{2}\right\rfloor$ identical vertices for the positive and negative literals respectively when the variable is universally quantified. Second, we need to avoid having the first burning vertex be fixed ahead of time. We thus take a large graph with a unique maximum length cooling sequence and replace the initial burning vertex with this graph. Finally, we need to take a winning strategy for the Saboteur and show that $\phi$ is satisfied. Our graph may satisfy this, but to make the proof more straightforward, we add a gadget that forces either $x_i$ or $\lnot x_i$ to be selected as the source in round $T+1 \leq i \leq T+n$, for some $T$ which will be determined later.

We will also briefly review the \emph{strong product} of graphs. The strong product of two graphs $G,H$ is a graph with vertex set $V(G) \times V(H)$ and $(u,x)$ adjacent to $(v,y)$ if $d_G(u,v)$ and $d_H(x,y)$ are both at most $1$. In the case where $H$ is a clique, this can be thought of as replacing each vertex of $G$ with a copy of $H$ and having all possible edges between cliques representing vertices that were adjacent in $G$.

\begin{proof}[Proof of Theorem~\ref{thm:pspace}]
Note that $k$-{\sc Liminal-Burning} is in PSPACE since the game lasts at most $|V(G)|$ rounds, so we need only show that it is PSPACE-hard. Let $Q_1 x_1 Q_2 x_2 \cdots Q_n x_n \phi(x_1,x_2,\dots,x_n)$ be an instance of {\sc $3$-QBF}, let $\phi$ have $m$ clauses, and let $G_{\phi}$ be as described above. 

We modify $G_\phi$ in the following way. For each $x_i$, add $k$-cliques, one for $x_i$ and one for $\lnot x_i$, then join them with all possible edges. Call these pairs of cliques \emph{double assignment gadgets}. For any literal $x_i$ and path from $x_i$ to $C_j$, add a path from the vertex adjacent to $C_j$ to the clique corresponding to $x_i$ with $m+2i-j-1$ internal vertices, and have the last vertex of the path, which is not in the clique, be adjacent to all vertices in the clique. We can assume every literal appears in a clause since adding clauses of the form $x_i \lor \lnot x_i \lor x_i$ never changes the truth value of $\phi$. For each variable $x_i,$ if $Q_i = \forall$ then replace the vertices corresponding to $x_i$ and $\lnot x_i$, respecting adjacency, with cliques of order $\left\lceil \frac{k}{2}\right\rceil$ and $\left\lfloor \frac{k}{2}\right\rfloor$ respectively. If $Q_i = \exists$, then replace the vertices corresponding to $x_i$ and $\lnot x_i$ with $k$-cliques. Additionally, replace each vertex in the clause gadgets with a $k$-clique. 

To complete our construction, let $T = \mathrm{diam}(G_\phi)$ and let $H_t$ be defined with $V(H_t) = \{u^i_j : 1\leq i \leq t, 1 \leq i \leq j\}$ and edges $u^{i}_a u^j_b$, where $i=j, \text { or } i = j+1 \text{ and } b < j, \text{ or } i+1 = j \text{ and } a < i.$ 
Cooling $u^i_i$ in round $i$ is the only optimal cooling strategy up to symmetry, so we can force the game to start with a certain sequence of moves. Now, delete $a_0$ from $G_\phi$, add in a copy of $H_T \boxtimes K_k$ and add an edge between $a_1$ and every $u^T_i$. Call the resulting graph $G'$, noting that $\mathrm{diam}(G')=T + 2n + m$. We claim that the instance of {\sc $3$-QBF} is true if and only if $b_k(G') \geq T + 2n + m + 1$. 

        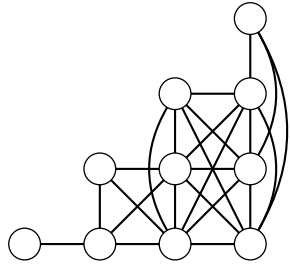
\begin{figure}[t]
            \centering 
            \begin{tikzpicture}
                \GraphInit[vstyle=Classic]
                \SetUpVertex[FillColor=white]
            
                \tikzset{VertexStyle/.append style={minimum size=12pt, inner sep=1pt}}
            
                \Vertex[x=0,y=0,NoLabel=true,]{V11}

                \Vertex[x=1,y=0,NoLabel=true,]{V21}
                \Vertex[x=1,y=1,NoLabel=true,]{V22}

                \Vertex[x=2,y=0,NoLabel=true,]{V31}
                \Vertex[x=2,y=1,NoLabel=true,]{V32}
                \Vertex[x=2,y=2,NoLabel=true,]{V33}

                \Vertex[x=3,y=0,NoLabel=true,]{V41}
                \Vertex[x=3,y=1,NoLabel=true,]{V42}
                \Vertex[x=3,y=2,NoLabel=true,]{V43}
                \Vertex[x=3,y=3,NoLabel=true,]{V44}

                \Edge(V22)(V21)

                \Edge(V31)(V32)
                \Edge(V32)(V33)

                \Edge(V41)(V42)
                \Edge(V42)(V43)
                \Edge(V43)(V44)
    
                \Edge(V11)(V21)

                \Edge(V21)(V31)
                \Edge(V21)(V32)
                \Edge(V22)(V31)
                \Edge(V22)(V32)

                \Edge(V31)(V41)
                \Edge(V31)(V42)
                \Edge(V31)(V43)

                \Edge(V32)(V41)
                \Edge(V32)(V42)
                \Edge(V32)(V43)

                \Edge(V33)(V41)
                \Edge(V33)(V42)
                \Edge(V33)(V43)

                \tikzset{EdgeStyle/.style = {bend left=30}}
                
                \Edge(V31)(V33)
                \Edge(V43)(V41)
                \Edge(V44)(V41)
                \Edge(V44)(V42)
                
            \end{tikzpicture}

            \caption{The graph $H_4$.}
        \end{figure}
        
If the {$3$-QBF} instance is true, then the Saboteur can make the game last $T + 2n + m + 1$ rounds. In round $i,$ for $1 \leq i \leq T$, the Saboteur reveals the clique corresponding to $u^i_i$. In round $T+i,$ for $1 \leq i \leq n$, the Saboteur has two strategies. If $Q_i = \exists$, then the Saboteur reveals the vertices in the clique corresponding to whichever of $x_i$ and $\lnot x_i$ is set to false. If $Q_i = \forall$, then the Saboteur reveals the vertices corresponding to $x_i$ and $\lnot x_i$, and whichever the Arsonist chooses, the formula will still be true. For the next $m$ rounds, the Saboteur reveals the clique corresponding to the true literal in each clause. Finally, for the final $n$ rounds, the Saboteur reveals the vertices in the clique corresponding to the value not taken by $x_i$, causing that clique not to burn yet. In round $T+2n + m + 1$, burning spreads to cover the final clique.
        
Suppose the game lasts $T + 2n + m + 1$ rounds. Since the diameter of $G_{\phi}$ is only witnessed by taking one vertex from the double assignment gadget for $x_n$ and one vertex from the clique corresponding to $u^1_1$, the first set revealed by the Saboteur must be a clique in either of these sets. Analogous to cooling, we can reverse the order of the sets revealed by the Saboteur here and still obtain an optimal strategy. Hence, we assume that the clique corresponding to $u^1_1$ is the first set revealed. Thus, the first $T$ rounds must consist of the Saboteur revealing the clique corresponding to $u^i_i$ in round $i$ for $1 \leq i \leq T$.  
        
In the final round, the spread phase covers the graph, so consider the $n$ rounds that precede this. It is evident that in round $T + n + m + j$ for $1 \leq i \leq n$, the Saboteur must reveal the cliques which prevent double assignment for $x_i$. Similarly, in round $T + n + j$ for $1 \leq j \leq m$, the Saboteur must reveal a clique in the gadget for $C_j$.

We derive the strategy for rounds $T + 1$ through $T + n$ so that these final $n + m$ choices are valid. First, suppose the Saboteur reveals a vertex on the path between a clause $C_j$ and a variable $x_i$. If this happens on or before round $T + i$, then $C_j$ is burned at the start of round $T + n + j$, so this cannot happen. Suppose this happens in round $T + k$ with $i < k \leq n$. If the game still lasts $T + 2n + m + 1$ rounds, then the variable vertex for $x_i$ was not chosen as a source previously, or else the gadget for $C_j$ would burn early. This ``sets'' $x_i$ to false in $C_j$, which causes the $x_i$ double assignment gadget to burn early since there are instances where $x_i$ is set to both false and true. Thus, the Saboteur cannot reveal any such vertices.

If the Saboteur reveals some $a_j$ with $j > i$ in round $T + i$, then the Saboteur must reveal at least one other vertex. If they are all from the $x_i$ variable gadget and burning $a_j$ is better for the Saboteur, then the Arsonist will burn one of the $x_i$ vertices. The Saboteur must also reveal some other $a_j$ or reveal a later variable vertex. In both cases, some clause will burn early. If the Saboteur reveals some vertex corresponding to $z_j$ with $j > i$ in round $T + i$, then the Arsonist will burn this vertex, and some clause will burn early. The Saboteur will only reveal vertices corresponding to $x_i$ or $\lnot x_i$ in round $T + i$ for $1 \leq i \leq n$. 

Finally, revealing some vertex on the paths from the clause gadgets to the double assignment gadgets will cause the double assignment gadget to burn early, which will make the Saboteur's strategy invalid. In an optimal strategy, when $Q_i = \exists$, the Saboteur will choose a value for $x_i$, and when $Q_i = \forall$, the Saboteur must let the Arsonist choose the value of $x_i$. Therefore, if the game lasts $T + 2n + m  + 1$ rounds, then all the clauses must be satisfied, and so the instance of {\sc $3$-QBF} must be true.
\end{proof}

If when we construct the graph $G'$ above we let $k$ be some large polynomial in the cardinality of the {\sc $3$-QBF} instance, then we can force $k$ (as a function of the number of vertices in the new graph) to be $n^\beta$ for any $\beta \in [0,1)$. 

\begin{corollary}
For all $\beta \in [0,1)$, $\left\lceil n^\beta \right\rceil$-{\sc Liminal-Burning} is PSPACE-complete.
\end{corollary}

We now consider the case $f(n) = \Omega(n)$. Define $K^t_r$ to be $t$ disjoint copies of $K_r$.
    
\begin{theorem} \label{thm:disj_cliques}
If $k \geq |V(G)|,$ then $b_k(G \cup K^t_k) = t + b(G)$.
\end{theorem}

\begin{proof}
We show that the Arsonist can always end the game in $t + b(G)$ rounds. Fix a vertex $u \in V(G)$ such that $u$ is an optimal first choice for burning on $G$. We claim that in the first $t + 1$ rounds, the Arsonist can always either burn a vertex in an unburned copy of $K_k$ or burn $u$. In round $1$, the Saboteur can either reveal all of $V(G)$ or reveal at least one vertex in one of the cliques, so we now proceed by induction. Suppose $u$ has already been burned in some previous round and that there are still unburned cliques. The Saboteur cannot reveal any vertices in a clique containing a source, and since $u$ has been burned, they can reveal at most $n-1 \leq k-1$ vertices in $V(G)$ and must reveal at least one vertex in an unburned clique. If $u$ is not yet burned, then the Saboteur can either reveal $u$ or reveal at most $n-1$ vertices in $V(G)$ and then must reveal a vertex in an unburned clique. 

Thus, at the end of round $t + 1$, every clique contains a source, and $u$ is a source. If $u$ was selected before round $t + 1$, then this aids the Arsonist, and by following a burning strategy on $G$, the Arsonist finishes burning the graph on or before round $t + b(G).$ If the Saboteur reveals only clique vertices until round $t + 1$ when they can only reveal vertices from $V(G)$, then this takes at least $t + b(G)$ rounds. 
\end{proof}

We next produce a reduction from {\sc Burning} to $k$-\textsc{Liminal-Burning} such that the instance of {\sc Burning} is true if and only if the instance of $k$-\textsc{Liminal-Burning} is false. This will show that $k$-\textsc{Liminal-Burning} is co-NP-hard.

\begin{theorem}
If $f(n) = \left\lceil \frac{n}{t} \right\rceil$ for some $t \in \mathbb{Z}^+$, then $f(n)$-\textsc{Liminal-Burning} is co-NP-hard.
\end{theorem}

\begin{proof}
Let $(G,c)$ be an instance of {\sc Burning}, let $n = |V(G)|$, and consider the instance $(G \cup K^{t-1}_{n},(t-1)+(c-1))$ of $f(n)$-\textsc{Liminal-Burning}. We then have that $f(|V(G \cup K^{t-1}_{n})|) = n$ and the result follows from Theorem~\ref{thm:disj_cliques}.
\end{proof}

\section{Future directions}
We introduced liminal burning, a two-player game on graphs that generalizes both burning and cooling. Determining $b_k$ for all values of $k$ for several families of graphs, such as hypercubes and paths, remains open. In particular, deriving tighter bounds for $b_k(Q_n)$ remains open. One direction is to study the game for random graphs, such as binomial random graphs or random regular graphs. Isoperimetry was an effective tool when studying the cooling number of grids in \cite{cooling}, and it may be useful for liminal burning in other graph families.

It would be interesting to exchange the goals of the Saboteur and the Arsonist, where the player who reveals vertices wants to minimize the game time, and the player who selects vertices wants to maximize it. We may also consider a non-monotonic version of the game as introduced in \cite{holden}, where revealed vertices become unrevealed in the next round. The complexity of $k$-\textsc{Liminal-Burning} remains open for certain orders of $k$, such as $k=\frac{n}{\log n}$ or more generally, if $k = 
\left(\bigcap_{\beta \in [0,1)}\omega\left(n^\beta\right)\right) \cap o(n).$

\end{document}